# RANK-BASED INFERENCE FOR BIVARIATE EXTREME-VALUE COPULAS

By Christian Genest[1] and Johan Segers[2]

*Université Laval and Université catholique de Louvain and Tilburg University*

Consider a continuous random pair $(X, Y)$ whose dependence is characterized by an extreme-value copula with Pickands dependence function $A$. When the marginal distributions of $X$ and $Y$ are known, several consistent estimators of $A$ are available. Most of them are variants of the estimators due to Pickands [*Bull. Inst. Internat. Statist.* **49** (1981) 859–878] and Capéraà, Fougères and Genest [*Biometrika* **84** (1997) 567–577]. In this paper, rank-based versions of these estimators are proposed for the more common case where the margins of $X$ and $Y$ are unknown. Results on the limit behavior of a class of weighted bivariate empirical processes are used to show the consistency and asymptotic normality of these rank-based estimators. Their finite- and large-sample performance is then compared to that of their known-margin analogues, as well as with endpoint-corrected versions thereof. Explicit formulas and consistent estimates for their asymptotic variances are also given.

**1. Introduction.** Let $(X, Y)$ be a pair of continuous random variables with joint and marginal distribution functions defined for all $x, y \in \mathbb{R}$ by

$$H(x, y) = \mathrm{P}(X \leq x, Y \leq y), \qquad F(x) = \mathrm{P}(X \leq x),$$
$$G(y) = \mathrm{P}(Y \leq y),$$

respectively. Let also $U = F(X)$ and $V = G(Y)$, and for all $u, v \in \mathbb{R}$, write $C(u, v) = \mathrm{P}(U \leq u, V \leq v)$. As shown by Sklar (1959), $C$ is the unique copula

Received June 2008; revised November 2008.

[1]Supported by grants from the Natural Sciences and Engineering Research Council of Canada, the Fonds québécois de la recherche sur la nature et les technologies and the Institut de finance mathématique de Montréal.

[2]Supported by the IAP research network Grant P6/03 of the Belgian Government (Belgian Science Policy).

for which $H$ admits the representation

$$(1.2) \qquad H(x, y) = C\{F(x), G(y)\}$$

for all $x, y \in \mathbb{R}$. The dependence between $X$ and $Y$ is characterized by $C$.

This paper is concerned with the estimation of $C$ under the assumption that it is an extreme-value copula, that is, when there exists a function $A : [0, 1] \to [1/2, 1]$ such that for all $(u, v) \in (0, 1)^2$,

$$(1.3) \qquad C(u, v) = (uv)^{A\{\log(v)/\log(uv)\}}.$$

It was shown by Pickands ([1981](#)) that $C$ is a copula if and only if $A$ is convex and $\max(t, 1 - t) \le A(t) \le 1$ for all $t \in [0, 1]$. By reference to this work, the function $A$ is often referred to as the "Pickands dependence function."

The interest in extreme-value copulas stems from their characterization as the large-sample limits of copulas of componentwise maxima of strongly mixing stationary sequences [Deheuvels ([1984](#)) and Hsing ([1989](#))]. More generally, these copulas provide flexible models for dependence between positively associated variables [Cebrián, Denuit and Lambert ([2003](#)), Ghoudi, Khoudraji and Rivest ([1998](#)) and Tawn ([1988](#))].

Parametric and nonparametric estimation methods for $A$ are reviewed in Section 9.3 of Beirlant et al. ([2004](#)). Nonparametric estimation based on a random sample $(X_1, Y_1), \ldots, (X_n, Y_n)$ from $H$ has been considered successively by Pickands ([1981](#)), Deheuvels ([1991](#)), Capéraà, Fougères and Genest ([1997](#)), Jiménez, Villa-Diharce and Flores ([2001](#)), Hall and Tajvidi ([2000](#)) and Segers ([2007](#)). In these papers, the margins $F$ and $G$ are assumed to be known, so that in effect, a random sample $(F(X_1), G(Y_1)), \ldots, (F(X_n), G(Y_n))$ from $C$ is available. In their multivariate extension of the Capéraà–Fougères–Genest (CFG) estimator, Zhang, Wells and Peng ([2008](#)) also assume knowledge of the marginal distributions.

In practice, however, margins are rarely known. A natural way to proceed is then to estimate $F$ and $G$ by their empirical counterparts $F_n$ and $G_n$, and to base the estimation of $C$ on the pseudo-observations $(F_n(X_1), G_n(Y_1)), \ldots, (F_n(X_n), G_n(Y_n))$. This amounts to working with the pairs of scaled ranks. This avenue was recently considered by Abdous and Ghoudi ([2005](#)), but only from a computational point of view. No theory was provided.

Rank-based versions of the Pickands and CFG estimators of $A$ are defined in Section [2](#). Endpoint corrections in the manner of Deheuvels and Hall–Tajvidi are also considered, but in contrast with the case of known margins, they have no effect asymptotically. Weak consistency and asymptotic normality of the rank-based estimators are established in Section [3](#). Explicit formulas and consistent estimates for their asymptotic variances are also given. Large-sample comparisons reported in Section [4](#) show that the rank-based estimators are more efficient than the uncorrected versions based on the true, known margins. Extensive numerical work also suggests that the CFG estimator is generally preferable to the Pickands estimator. A small simulation study reported in Section 5 provides evidence that the conclusions remain valid in small samples, and a few closing comments are made in Section 6.

In order to ease reading, all technical arguments are relegated to Appendices A–F. The developments rely heavily on a limit theorem for weighted bivariate empirical processes, which may be of independent interest. A statement and proof of the latter result are given in Appendix G.

The following notational conventions are used in the sequel. For $x, y \in \mathbb{R}$, let $x \wedge y = \min(x, y)$ and $x \vee y = \max(x, y)$. The arrow $\rightsquigarrow$ means weak convergence, and $\mathbf{1}(E)$ stands for the indicator function of the event $E$. Given a univariate cumulative distribution function $F$, its left-continuous generalized inverse is denoted by $F^{\leftarrow}$. Furthermore, $\mathcal{C}([0, 1])$ represents the space of continuous, real-valued functions on $[0, 1]$, while $\ell^{\infty}(\mathcal{W})$ is the space of bounded, real-valued functions on the set $\mathcal{W}$; both are equipped with the topology of uniform convergence.

## 2. Estimators of the dependence function.

Consider a pair $(X, Y)$ of continuous random variables whose joint distribution function $H$ has margins $F$ and $G$, as per (1.1). Assume that the unique function $C$ implicitly defined by (1.2) belongs to the class (1.3) of extreme-value copulas.

Let $U = F(X)$ and $V = G(Y)$. The pair $(U, V)$ is then distributed as $C$. Accordingly, the variables $S = -\log U$ and $T = -\log V$ are exponential with unit mean. For all $t \in (0, 1)$, write

$$\xi(t) = \frac{S}{1-t} \wedge \frac{T}{t}$$

and set $\xi(0) = S$, $\xi(1) = T$. Note that for any $t \in [0, 1]$ and $x \geq 0$, one has

$$P\{\xi(t) > x\} = P\{S > (1-t)x, T > tx\}$$
$$= P\{U < e^{-(1-t)x}, V < e^{-tx}\} = e^{-xA(t)}.$$

Thus, $\xi(t)$ is an exponential random variable with

$$(2.1) \qquad E\{\xi(t)\} = 1/A(t) \quad \text{and} \quad E\{\log \xi(t)\} = -\log A(t) - \gamma,$$

where $\gamma = -\int_0^{\infty} \log(x) e^{-x} \, dx \approx 0.577$ is Euler's constant. These observations motivate the following estimators of $A$.

### 2.1. Pickands and CFG estimators.

Let $(X_1, Y_1), \ldots, (X_n, Y_n)$ be a random sample from $H$. For $i \in \{1, \ldots, n\}$, let $U_i = F(X_i)$, $V_i = G(Y_i)$, $S_i = -\log(U_i) = \xi_i(0)$ and $T_i = -\log(V_i) = \xi_i(1)$, and for all $t \in (0, 1)$, write

$$\xi_i(t) = \frac{S_i}{1-t} \wedge \frac{T_i}{t}.$$

When the margins $F$ and $G$ are known, the estimation of $A(t)$ for arbitrary $t \in (0,1)$ can be based on the sample $\xi_1(t), \ldots, \xi_n(t)$. In view of (2.1), two obvious solutions are defined implicitly by

$$1/A_{n,\mathrm{u}}^{\mathrm{P}}(t) = \frac{1}{n} \sum_{i=1}^{n} \xi_i(t),$$

$$\log A_{n,\mathrm{u}}^{\mathrm{CFG}}(t) = -\gamma - \frac{1}{n} \sum_{i=1}^{n} \log \xi_i(t).$$

The first estimator was proposed by Pickands (1981). As shown by Segers (2007), the second is an uncorrected form of the estimator due to Capéraà, Fougères and Genest (1997). The index "u" is used to distinguish these estimators from endpoint-corrected versions introduced in Section 2.3.

2.2. *Rank-based versions of the Pickands and CFG estimators.* When the margins $F$ and $G$ are unknown, a natural solution is to rely on their empirical counterparts, $F_n$ and $G_n$. This leads to pairs $(F_n(X_1), G_n(Y_1)), \ldots,$ $(F_n(X_n), G_n(Y_n))$ of pseudo-observations for $C$. To avoid dealing with points at the boundary of the unit square, however, it is more convenient to work with scaled variables $\hat{U}_i = nF_n(X_i)/(n+1)$ and $\hat{V}_i = nG_n(Y_i)/(n+1)$ defined explicitly for every $i \in \{1, \ldots, n\}$ by

$$\hat{U}_i = \frac{1}{n+1} \sum_{j=1}^{n} \mathbf{1}(X_j \leq X_i), \qquad \hat{V}_i = \frac{1}{n+1} \sum_{j=1}^{n} \mathbf{1}(Y_j \leq Y_i).$$

The pair $(\hat{U}_i, \hat{V}_i)$, whose coordinates are scaled ranks, can be regarded as a sample analogue of the unobservable pair $(U_i, V_i) = (F(X_i), G(Y_i))$. For every $i \in \{1, \ldots, n\}$ and for arbitrary $t \in (0,1)$, let

$$\hat{S}_i = -\log \hat{U}_i = \hat{\xi}_i(0), \qquad \hat{T}_i = -\log \hat{V}_i = \hat{\xi}_i(1), \qquad \hat{\xi}_i(t) = \frac{\hat{S}_i}{1-t} \wedge \frac{\hat{T}_i}{t}.$$

Rank-based versions of $A_{n,\mathrm{u}}^{\mathrm{P}}(t)$ and $A_{n,\mathrm{u}}^{\mathrm{CFG}}(t)$ are then given by

$$(2.2) \qquad 1/A_{n,\mathrm{r}}^{\mathrm{P}}(t) = \frac{1}{n} \sum_{i=1}^{n} \hat{\xi}_i(t),$$

$$(2.3) \qquad \log A_{n,\mathrm{r}}^{\mathrm{CFG}}(t) = -\gamma - \frac{1}{n} \sum_{i=1}^{n} \log \hat{\xi}_i(t).$$

It is these estimators that are the focus of the present study. Before proceeding, however, endpoint corrections will be discussed briefly.

2.3. *Endpoint corrections to the Pickands estimator.* The Pickands estimator does not satisfy the endpoint constraints

$$(2.4) \qquad\qquad A(0) = A(1) = 1.$$

In order to overcome this defect, Deheuvels ([1991](#)) introduced the estimator

$$1/A_{n,c}^{P}(t) = 1/A_{n,u}^{P}(t) - (1-t)\{1/A_{n,u}^{P}(0) - 1\} - t\{1/A_{n,u}^{P}(1) - 1\}.$$

More generally, Segers ([2007](#)) considered endpoint corrections of the form

$$1/A_{n,ab}^{P} = 1/A_{n,u}^{P}(t) - a(t)\{1/A_{n,u}^{P}(0) - 1\} - b(t)\{1/A_{n,u}^{P}(1) - 1\},$$

where $a, b : [0, 1] \to \mathbb{R}$ are arbitrary continuous mappings. He showed how to select (and estimate) $a$ and $b$ so as to minimize the asymptotic variance of the endpoint-corrected estimator at every $t \in [0, 1]$. In particular, he showed that $a(t) = 1 - t$, $b(t) = t$ are optimal at independence. For arbitrary $A$, the Pickands estimator with optimal correction is denoted $A_{n,\text{opt}}^{P}$ hereafter.

Similar strategies can be used to ensure that the rank-based estimator $A_{n,r}^{P}$ fulfills conditions ([2.4](#)). Note, however, that whatever the choice of mappings $a$ and $b$, the endpoint-corrected estimator is then asymptotically equivalent to $A_{n,r}^{P}$, because

$$(2.5) \quad 1/A_{n,r}^{P}(0) = 1/A_{n,r}^{P}(1) = \frac{1}{n}\sum_{i=1}^{n}\log\{(n+1)/i\} = 1 + O(n^{-1}\log n).$$

An alternative correction inspired by Hall and Tajvidi ([2000](#)) is given by

$$1/A_{n,r}^{\text{HT}}(t) = \frac{1}{n}\sum_{i=1}^{n}\bar{\xi}_i(t),$$

where for all $t \in [0, 1]$ and $i \in \{1, \ldots, n\}$,

$$\bar{\xi}_i(t) = \frac{\bar{S}_i}{1-t} \wedge \frac{\bar{T}_i}{t}$$

with $\bar{S}_i = n\hat{S}_i/(\hat{S}_1 + \cdots + \hat{S}_n)$ and $\bar{T}_i = n\hat{T}_i/(\hat{T}_1 + \cdots + \hat{T}_n)$.

By construction, one has $A_{n,r}^{\text{HT}}(0) = A_{n,r}^{\text{HT}}(1) = 1$, but an additional merit of this estimator is that $A_{n,r}^{\text{HT}}(t) \geq t \vee (1-t)$ for all $t \in [0, 1]$. Note, however, that because $A_{n,r}^{\text{HT}} = A_{n,r}^{P}/A_{n,r}^{P}(0)$, this estimator is again asymptotically indistinguishable from $A_{n,r}^{P}$ in view of ([2.5](#)).

2.4. *Endpoint corrections to the CFG estimator.* In order to meet constraints ([2.4](#)), Capéraà, Fougères and Genest ([1997](#)) consider estimators of the form

$$\log A_{n,c}^{\text{CFG}}(t) = \log A_{n,u}^{\text{CFG}}(t) - p(t)\log A_{n,u}^{\text{CFG}}(0) - \{1 - p(t)\}\log A_{n,u}^{\text{CFG}}(1),$$

where $p\colon [0,1] \to \mathbb{R}$ is an arbitrary continuous mapping. They use $p(t) = 1 - t$ as an expedient when the margins are known. The more general choice,

$$\log A_{n,\mathrm{ab}}^{\mathrm{CFG}}(t) = \log A_{n,\mathrm{u}}^{\mathrm{CFG}}(t) - a(t) \log A_{n,\mathrm{u}}^{\mathrm{CFG}}(0) - b(t) \log A_{n,\mathrm{u}}^{\mathrm{CFG}}(1)$$

is investigated by Segers ([2007](#)), who identified the optimal functions $a$ and $b$. The resulting estimator is hereafter denoted $A_{n,\mathrm{opt}}^{\mathrm{CFG}}$.

When the margins are unknown, however, the correction to the estimator $A_{n,\mathrm{r}}^{\mathrm{CFG}}$ (and hence the choices of $p$, $a$ and $b$) has no impact on the asymptotic distribution of this rank-based statistic, because

$$-\log A_{n,\mathrm{r}}^{\mathrm{CFG}}(0) = -\log A_{n,\mathrm{r}}^{\mathrm{CFG}}(1)$$
$$= \frac{1}{n} \sum_{i=1}^{n} \log\log\{(n+1)/i\} - \int_0^1 \log\log(1/x)\,dx$$
$$= O\{n^{-1}(\log n)^2\}.$$

**3. Asymptotic results.** The limiting behavior of the estimators $A_{n,\mathrm{r}}^{\mathrm{P}}$ and $A_{n,\mathrm{r}}^{\mathrm{CFG}}$ (and of their asymptotically equivalent variants) can be determined once they have been expressed as appropriate functionals of the empirical copula $C_n$, defined for all $u, v \in [0,1]$ by

$$(3.1) \qquad \hat{C}_n(u,v) = \frac{1}{n} \sum_{i=1}^{n} \mathbf{1}(\hat{U}_i \leq u, \hat{V}_i \leq v).$$

Note that this definition is somewhat different from the original one given by Deheuvels ([1979](#)). From inequality ([B.1](#)) in Appendix [B](#), one can see that the difference between the two versions is $O(n^{-1})$ as $n \to \infty$ almost surely.

The following lemma is proved in Appendix [A](#).

LEMMA 3.1. *For every $t \in [0,1]$, one has*

$$(3.2) \qquad 1/A_{n,\mathrm{r}}^{\mathrm{P}}(t) = \int_0^1 \hat{C}_n(u^{1-t}, u^t)\,\frac{du}{u},$$

$$(3.3) \qquad \log A_{n,\mathrm{r}}^{\mathrm{CFG}}(t) = -\gamma + \int_0^1 \{\hat{C}_n(u^{1-t}, u^t) - \mathbf{1}(u > e^{-1})\}\,\frac{du}{u \log u}.$$

Note that replacing $\hat{C}_n(u^{1-t}, u^t)$ by $C(u^{1-t}, u^t) = u^{A(t)}$ in ([3.2](#)) and ([3.3](#)) yields $1/A(t)$ and $\log A(t)$, respectively. Therefore, both $A_{n,\mathrm{r}}^{\mathrm{P}}$ and $A_{n,\mathrm{r}}^{\mathrm{CFG}}$ can be expected to yield consistent and asymptotically unbiased estimators of $A$. More generally, the asymptotic behavior of the processes

$$\mathbb{A}_{n,\mathrm{r}}^{\mathrm{P}} = n^{1/2}(A_{n,\mathrm{r}}^{\mathrm{P}} - A) \quad \text{and} \quad \mathbb{A}_{n,\mathrm{r}}^{\mathrm{CFG}} = n^{1/2}(A_{n,\mathrm{r}}^{\mathrm{CFG}} - A)$$

is a function of the limit, $\mathbb{C}$, of the empirical copula process

$$(3.4) \qquad \mathbb{C}_n = n^{1/2}(\hat{C}_n - C).$$

3.1. *Limiting behavior of $A_{n,r}^P$ and $A_{n,r}^{CFG}$.* As shown under various conditions by Rüschendorf ([1976](#)), Stute ([1984](#)), Fermanian, Radulović and Wegkamp ([2004](#)) and Tsukahara ([2005](#)), the weak limit $\mathbb{C}$ of the process $\mathbb{C}_n$ is closely related to a bivariate pinned $C$-Brownian sheet $\alpha$, that is, a centered Gaussian random field on $[0,1]^2$ whose covariance function is defined for every value of $u, v, u', v' \in [0,1]$ by

$$\mathrm{cov}\{\alpha(u,v), \alpha(u',v')\} = C(u \wedge u', v \wedge v') - C(u,v)C(u',v').$$

Denoting $\dot{C}_1(u,v) = \partial C(u,v)/\partial u$ and $\dot{C}_2(u,v) = \partial C(u,v)/\partial v$, one has

$$\mathbb{C}(u,v) = \alpha(u,v) - \dot{C}_1(u,v)\alpha(u,1) - \dot{C}_2(u,v)\alpha(1,v)$$

for every pair $(u,v) \in [0,1]^2$. The weak limits of the rank-based processes $A_{n,r}^P$ and $A_{n,r}^{CFG}$ are then, respectively, defined at each $t \in [0,1]$ by

$$(3.5) \qquad \mathbb{A}_r^P(t) = -A^2(t) \int_0^1 \mathbb{C}(u^{1-t}, u^t) \, \frac{du}{u},$$

$$(3.6) \qquad \mathbb{A}_r^{CFG}(t) = A(t) \int_0^1 \mathbb{C}(u^{1-t}, u^t) \, \frac{du}{u \log u}.$$

This fact, which is the main result of the present paper, is stated formally below under the assumption that $A$ is twice continuously differentiable. This hypothesis could possibly be relaxed, but at the cost of an extension of the strong approximation results in Stute ([1984](#)) and Tsukahara ([2005](#)).

THEOREM 3.2. *If $A$ is twice continuously differentiable, then $A_{n,r}^P \rightsquigarrow \mathbb{A}_r^P$ and $A_{n,r}^{CFG} \rightsquigarrow \mathbb{A}_r^{CFG}$ as $n \to \infty$ in the space $\mathcal{C}([0,1])$ equipped with the topology of uniform convergence.*

This result, which is proved in Appendix [B](#), is to be contrasted with the case of known margins, where one has access to the pairs $(U_i, V_i) = (F(X_i), G(Y_i))$ for all $i \in \{1, \dots, n\}$. As shown by Segers ([2007](#)), the estimators $A_{n,u}^P$ and $A_{n,u}^{CFG}$ are then of the same form as in ([3.2](#)) and ([3.3](#)), but with $\hat{C}_n$ replaced by the empirical distribution function

$$C_n(u,v) = \frac{1}{n} \sum_{i=1}^n \mathbf{1}(U_i \leq u, V_i \leq v).$$

The asymptotic behavior of the estimators is then as in Theorem [3.2](#), but with the process $\mathbb{C}$ in ([3.5](#)) and ([3.6](#)) replaced by the $C$-Brownian sheet $\alpha$. In what follows, the weak limits of the processes $n^{1/2}(A_{n,u}^P - A)$ and $n^{1/2}(A_{n,u}^{CFG} - A)$ are denoted $\mathbb{A}_u^P$ and $\mathbb{A}_u^{CFG}$, respectively.

3.2. *Asymptotic variances of $A_{n,r}^P$ and $A_{n,r}^{CFG}$.* Fix $u, v, t \in [0,1]$ and let

$$\sigma(u,v;t) = \mathrm{cov}\{\mathbb{C}(u^{1-t}, u^t), \mathbb{C}(v^{1-t}, v^t)\}.$$

In view of Theorem 3.2, the asymptotic variances of the estimators $A_{n,r}^P(t)$ and $A_{n,r}^{CFG}(t)$ of $A(t)$ are given by

$$\mathrm{var}\,\mathbb{A}_r^P(t) = A^4(t) \int_0^1 \int_0^1 \sigma(u,v;t)\,\frac{du}{u}\,\frac{dv}{v},$$

$$\mathrm{var}\,\mathbb{A}_r^{CFG}(t) = A^2(t) \int_0^1 \int_0^1 \sigma(u,v;t)\,\frac{du}{u \log u}\,\frac{dv}{v \log v}.$$

Closed-form expressions for the latter are given next in terms of

$$\mu(t) = A(t) - t\dot{A}(t), \qquad \nu(t) = A(t) + (1-t)\dot{A}(t),$$

where $\dot{A}(t) = dA(t)/dt \in [-1,1]$ for all $t \in (0,1)$. The asymptotic variance of $A_{n,r}^{CFG}$ also involves the dilogarithm function $L_2$, defined for $x \in [-1,1]$ by

$$L_2(x) = -\int_0^x \log(1-z)\,\frac{dz}{z} = \sum_{k=1}^\infty \frac{x^k}{k^2}.$$

The proof of the following result is presented in Appendix C.

PROPOSITION 3.3. *For $t \in [0,1]$, let $A_1(t) = A(t)/t$, $A_2(t) = A(t)/(1-t)$, $\bar{\mu}(t) = 1 - \mu(t)$ and $\bar{\nu}(t) = 1 - \nu(t)$. Then, $A^{-2}(t)\,\mathrm{var}\,\mathbb{A}_r^P(t)$ is given by*

$$2 - \{\mu(t) + \nu(t) - 1\}^2 - \frac{2\mu(t)\bar{\mu}(t)A_2(t)}{2A_2(t) - 1} - \frac{2\nu(t)\bar{\nu}(t)A_1(t)}{2A_1(t) - 1}$$

$$+ 2\mu(t)\nu(t)A_1(t)A_2(t)\int_0^1 \{A(s) + sA_1(t) + (1-s)A_2(t) - 1\}^{-2}\,ds$$

$$- 2\mu(t)A_1(t)A_2(t)\int_0^t [A(s) + (1-s)\{A_2(t) - 1\}]^{-2}\,ds$$

$$- 2\nu(t)A_1(t)A_2(t)\int_t^1 [A(s) + s\{A_1(t) - 1\}]^{-2}\,ds,$$

*while $A^{-2}(t)\,\mathrm{var}\,\mathbb{A}_r^{CFG}(t)$ equals*

$$\{1 + \mu^2(t) + \nu^2(t) - \mu(t) - \nu(t)\}L_2(1)$$

$$- 2\mu(t)\bar{\mu}(t)L_2\{-1 + 1/A_2(t)\} - 2\nu(t)\bar{\nu}(t)L_2\{-1 + 1/A_1(t)\}$$

$$- 2\mu(t)\nu(t)\int_0^1 \log\left\{1 - \frac{t(1-t)}{A(t)}\frac{1 - A(s)}{t(1-s) + (1-t)s}\right\}\frac{ds}{s(1-s)}$$

$$+ 2\mu(t)\int_0^t \log\left[1 - \frac{t(1-t)}{A(t)}\frac{1 - A(s) + s\{A_1(t) - 1\}}{t(1-s) + (1-t)s}\right]\frac{ds}{s(1-s)}$$

$$+ 2\nu(t)\int_t^1 \log\left[1 - \frac{t(1-t)}{A(t)}\frac{1 - A(s) + (1-s)\{A_2(t) - 1\}}{t(1-s) + (1-t)s}\right]\frac{ds}{s(1-s)}.$$

As stated below, great simplifications occur at independence. The proof of this result is given in Appendix D.

COROLLARY 3.4. *If $A \equiv 1$, then $\mu \equiv \nu \equiv 1$, and for all $t \in [0, 1]$,*

$$\operatorname{var} \mathbb{A}_r^P(t) = \frac{3t(1-t)}{(2-t)(1+t)},$$

$$\operatorname{var} \mathbb{A}_r^{CFG}(t) = 2L_2(-1) - 2L_2(t-1) - 2L_2(-t).$$

3.3. *Consistent estimates of the asymptotic variances.* Proposition 3.3 can be used to construct consistent estimates of $\operatorname{var} \mathbb{A}_r^P(t)$ and $\operatorname{var} \mathbb{A}_r^{CFG}(t)$ for arbitrary $t \in [0, 1]$. The latter is useful, for example, for the construction of asymptotic confidence intervals for the Pickands dependence function.

Specifically, suppose that $(A_n)$ is any sequence of consistent estimators for $A$, that is, suppose that $\|A_n - A\| \to 0$ in probability as $n \to \infty$, where $\|\cdot\|$ is the supremum norm on $\mathcal{C}([0, 1])$. Put

$$\hat{A}_n = \text{greatest convex minorant of } (A_n \wedge 1) \vee I \vee (1 - I),$$

where $I$ denotes the identity function. One can then invoke a lemma of Marshall (1970) to deduce that $\|\hat{A}_n - A\| \leq \|A_n - A\|$ and, hence, that $(\hat{A}_n)$ is also a consistent sequence of estimators. Furthermore, $\hat{A}_n$ is itself a Pickands dependence function. See Fils-Villetard, Guillou and Segers (2008) for another way of converting a pilot estimate $A_n$ into a Pickands dependence function.

Let $\hat{A}_n'$ denote the right-hand side derivative of $\hat{A}_n$. Because $\hat{A}_n$ is convex for every $n \in \mathbb{N}$, it is not hard to see that if $A$ is continuously differentiable in $t \in (0, 1)$, then $\hat{A}_n'(t) \to \dot{A}(t)$ in probability as $n \to \infty$. Consequently, if $A$ is replaced by $\hat{A}_n$ in Proposition 3.3, it can be seen that the resulting expressions converge in probability. In other words,

$$n \operatorname{var} \hat{A}_{n,r}^P(t) \to \operatorname{var} \mathbb{A}_r^P(t) \quad \text{and} \quad n \operatorname{var} A_{n,r}^{CFG}(t) \to \operatorname{var} \mathbb{A}_r^{CFG}(t).$$

4. **Efficiency comparisons.** Which of the rank-based estimators $A_{n,r}^P$ and $A_{n,r}^{CFG}$ is preferable in practice? When the margins are known, how do they fare compared with their uncorrected, corrected and optimal competitors? These issues are considered next in terms of asymptotic efficiency.

Figure 1 summarizes the findings, based either on mathematical derivations or on numerical calculations. In the diagram, an arrow $E_1 \to E_2$ between estimators $E_1$ and $E_2$ means that the latter is asymptotically more efficiency than the former, that is, $\sigma_{E_2}^2(t) \leq \sigma_{E_1}^2(t)$ for all $t \in [0, 1]$.

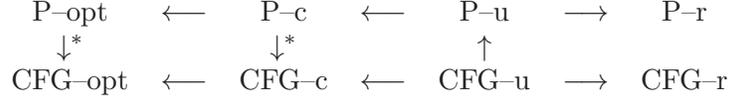

FIG. 1. *Comparisons of estimators. An arrow $E_1 \rightarrow E_2$ between estimators $E_1$ and $E_2$ means that the latter is asymptotically more efficiency than the former, that is, $\sigma_{E_2}^2 \leq \sigma_{E_1}^2$. Notation:* P = *Pickands;* CFG = *Capéraà–Fougères–Genest;* −u = *known margins, uncorrected;* −c = *known margins, endpoint-corrected;* −opt = *known margins, optimal correction;* −r = *rank-based. No arrow can be drawn between* P–r *and* CFG–r, *between* P–c *and* P–r, *or between* CFG–c *and* CFG–r. *An asterisk* (\*) *marks a conjecture based on extensive numerical computations.*

### 4.1. Uncorrected versus corrected estimators.

For all $u, v, t \in [0, 1]$, let

$$(4.1) \quad \sigma_0(u, v; t) = \mathrm{cov}\{\alpha(u^{1-t}, u^t), \alpha(v^{1-t}, v^t)\} = (u \wedge v)^{A(t)} - (uv)^{A(t)}.$$

It follows from the work of Segers (2007) that the asymptotic variances of the raw estimators $A_{n,u}^{\mathrm{P}}(t)$ and $A_{n,u}^{\mathrm{CFG}}(t)$ are given by

$$\mathrm{var}\,\mathbb{A}_u^{\mathrm{P}}(t) = A^4(t) \int_0^1 \int_0^1 \sigma_0(u, v; t)\, \frac{du}{u}\, \frac{dv}{v} = A^2(t),$$

$$\mathrm{var}\,\mathbb{A}_u^{\mathrm{CFG}}(t) = A^2(t) \int_0^1 \int_0^1 \sigma_0(u, v; t)\, \frac{du}{u \log u}\, \frac{dv}{v \log v} = L_2(1) A^2(t),$$

respectively. As $L_2(1) = \pi^2/6$, Pickands' original estimator $A_{n,u}^{\mathrm{P}}$ is more efficient than the uncorrected CFG estimator, $A_{n,u}^{\mathrm{CFG}}$, that is, CFG–u → P–u.

Formulas for the asymptotic variances of the endpoint-corrected versions $A_{n,c}^{\mathrm{P}}$ and $A_{n,c}^{\mathrm{CFG}}$ are more complex. They can be derived from the fact, also established by Segers (2007), that for all choices of continuous mappings $a, b : [0, 1] \rightarrow \mathbb{R}$, the weak limit $\mathbb{A}_{\mathrm{ab}}^{\mathrm{P}}$ of the process $n^{1/2}(A_{n,\mathrm{ab}}^{\mathrm{P}} - A)$ satisfies

$$(4.2) \qquad \mathbb{A}_{\mathrm{ab}}^{\mathrm{P}}(t) = \mathbb{A}_u^{\mathrm{P}}(t) - a(t)\mathbb{A}_u^{\mathrm{P}}(0) - b(t)\mathbb{A}_u^{\mathrm{P}}(1)$$

for all $t \in [0, 1]$. A similar result holds for the limit of $n^{1/2}(A_{n,\mathrm{opt}}^{\mathrm{CFG}} - A)$.

Using these facts, one can show that asymptotic variance reduction results from the application of the endpoint correction $a(t) = 1 - t$, $b(t) = t$, both for the Pickands and CFG estimators. In other words, one has P–u → P–c and CFG–u → CFG–c. This fact is formally stated below.

PROPOSITION 4.1. *For all choices of $A$ and $t \in [0, 1]$, one has*

$$\mathrm{var}\,\mathbb{A}_c^{\mathrm{P}}(t) \leq \mathrm{var}\,\mathbb{A}_u^{\mathrm{P}}(t) \quad and \quad \mathrm{var}\,\mathbb{A}_c^{\mathrm{CFG}}(t) \leq \mathrm{var}\,\mathbb{A}_u^{\mathrm{CFG}}(t).$$

The proof of this result may be found in Appendix E. Of course, it is trivial that P–c → P–opt and CFG–c → CFG–opt.

4.2. *Rank-based versus uncorrected estimators.* Somewhat more surprising, perhaps, is the fact that the rank-based versions of the Pickands and CFG estimators are asymptotically more efficient than their uncorrected, known-margin counterparts, that is, P–u → P–r and CFG–u → CFG–r. This observation is a consequence of the following more general result, whose proof may be found in Appendix F.

PROPOSITION 4.2. *For all choices of $A$ and $u, v, t \in [0, 1]$, one has*

$$\mathrm{cov}\{\mathbb{C}(u^{1-t}, u^t), \mathbb{C}(v^{1-t}, v^t)\} \leq \mathrm{cov}\{\alpha(u^{1-t}, u^t), \alpha(v^{1-t}, v^t)\}.$$

Consequently, the use of ranks improves the asymptotic efficiency of *any* estimator of $A$ whose limiting variance depends on $\sigma_0$ as defined in (4.1) through an expression of the form

$$\int_0^1 \int_0^1 \sigma_0(u, v; t) f(u, v) \, du \, dv,$$

where $f$ is nonnegative. See Henmi (2004) for other cases where efficiency is improved through the estimation of nuisance parameters whose value is known. Typically, this phenomenon occurs when the initial estimator is not semiparametrically efficient. Apparently, such is the case here, both for $A_{n,\mathrm{u}}^{\mathrm{P}}$ and $A_{n,\mathrm{u}}^{\mathrm{CFG}}$. As will be seen below, however, the phenomenon does not persist when endpoint-corrected estimators are considered.

4.3. *Ranked-based versus optimally corrected estimators.* Although efficiency comparisons between ranked-based and uncorrected versions of the Pickands and CFG estimators are interesting from a philosophical point of view, endpoint-corrected versions are preferable to the uncorrected estimators when margins are known. For this reason, comparisons between rank-based and corrected estimators are more relevant.

To investigate this issue, plots of the asymptotic variances of

$$A_{n,\mathrm{c}}^{\mathrm{P}}, \qquad A_{n,\mathrm{opt}}^{\mathrm{P}}, \qquad A_{n,\mathrm{r}}^{\mathrm{P}}, \qquad A_{n,\mathrm{c}}^{\mathrm{CFG}}, \qquad A_{n,\mathrm{opt}}^{\mathrm{CFG}}, \qquad A_{n,\mathrm{r}}^{\mathrm{CFG}}$$

were drawn for the following extreme-value copula models:

(a) The independence model, that is, $A(t) = 1$ for all $t \in [0, 1]$.
(b) The asymmetric logistic model [Tawn (1988)], namely,

$$A(t) = (1 - \psi_1)(1 - t) + (1 - \psi_2)t + [(\psi_1 t)^{1/\theta} + \{\psi_2(1 - t)\}^{1/\theta}]^\theta$$

with parameters $\theta \in (0, 1]$, $\psi_1, \psi_2 \in [0, 1]$. The special case $\psi_1 = \psi_2 = 1$ corresponds to the (symmetric) model of Gumbel (1960).

(c) The asymmetric negative logistic model [Joe (1990)], namely,

$$A(t) = 1 - [\{\psi_1(1-t)\}^{-1/\theta} + (\psi_2 t)^{-1/\theta}]^{-\theta}$$

with parameters $\theta \in (0, \infty)$, $\psi_1, \psi_2 \in (0, 1]$. The special case $\psi_1 = \psi_2 = 1$ gives the (symmetric) negative logistic of Galambos (1978).

(d) The asymmetric mixed model [Tawn (1988)], namely,

$$A(t) = 1 - (\theta + \kappa)t + \theta t^2 + \kappa t^3$$

with parameters $\theta$ and $\kappa$ satisfying $\theta \geq 0$, $\theta + 3\kappa \geq 0$, $\theta + \kappa \leq 1$, $\theta + 2\kappa \leq 1$. The special case $\kappa = 0$ and $\theta \in [0, 1]$ yields the (symmetric) mixed model [Tiago de Oliveira (1980)].

(e) The bilogistic model [Coles and Tawn (1994) and Joe, Smith and Weissman (1992)], namely,

$$A(t) = \int_0^1 \max\{(1-\beta)w^{-\beta}(1-t), (1-\delta)(1-w)^{-\delta}t\}\, dw$$

with parameters $(\beta, \delta) \in (0, 1)^2 \cup (-\infty, 0)^2$.

(f) The model of Hüsler and Reiss (1989), namely,

$$A(t) = (1-t)\Phi\left(\lambda + \frac{1}{2\lambda}\log\frac{1-t}{t}\right) + t\Phi\left(\lambda + \frac{1}{2\lambda}\log\frac{t}{1-t}\right),$$

where $\lambda \in (0, \infty)$ and $\Phi$ is the standard normal distribution function.

(g) The $t$-EV model [Demarta and McNeil (2005)], in which

$$A(w) = wt_{\chi+1}(z_w) + (1-w)t_{\chi+1}(z_{1-w}),$$

$$z_w = (1+\chi)^{1/2}[\{w/(1-w)\}^{1/\chi} - \rho](1-\rho^2)^{-1/2}$$

with parameters $\chi > 0$ and $\rho \in (-1, 1)$, where $t_{\chi+1}$ is the distribution function of a Student-$t$ random variable with $\chi + 1$ degrees of freedom.

Figure 2 corresponds to the case of independence. One can see from it that when $A \equiv 1$, the rank-based versions of the Pickands and CFG estimators are more efficient than the corresponding optimal, endpoint-corrected versions, even though the latter use information about the margins. As illustrated in Figure 3, however, the rank-based estimators are not always superior. The paradoxical phenomenon mentioned in Section 4.2 thus vanishes. Further, note the following:

(a) When $A$ is symmetric, one would expect the asymptotic variance of an estimator to reach its maximum at $t = 1/2$. Such is not always the case, however, as illustrated by the $t$-EV model.

(b) In the asymmetric negative logistic model, the asymptotic variance of the rank-based and optimally endpoint-corrected estimators is close to zero for all $t \in [0, 0.3]$. This is due to the fact that $A(t) \approx 1 - t$ on this interval when $\theta = 1/10$, $\psi_1 = 1/2$ and $\psi_2 = 1$ in this model.

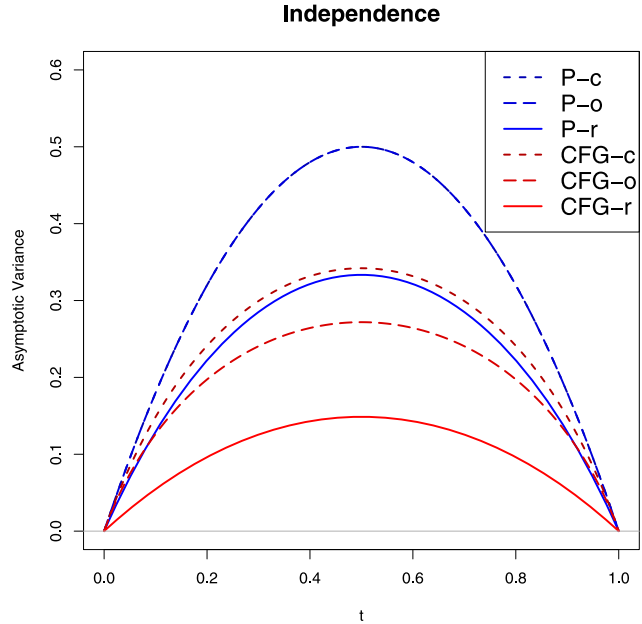

Fig. 2. *Graph, as a function of $t$, of the asymptotic variances of the estimators $A_{n,c}^{\mathrm{P}}(t) = A_{n,\mathrm{opt}}^{\mathrm{P}}(t)$, $A_{n,\mathrm{r}}^{\mathrm{P}}(t)$, $A_{n,c}^{\mathrm{CFG}}(t)$, $A_{n,\mathrm{opt}}^{\mathrm{CFG}}(t)$ and $A_{n,\mathrm{r}}^{\mathrm{CFG}}(t)$ in case of independence, $A \equiv 1$.*

4.4. *Comparison between the Pickands and CFG estimators.* In view of Figure 3, neither of the two rank-based estimators dominates the other one in terms of asymptotic efficiency. In most cases, however, the CFG estimator is superior to the Pickands estimator. In the figure, $A_{n,\mathrm{opt}}^{\mathrm{CFG}}$ is also seen to be systematically more efficient than $A_{n,\mathrm{opt}}^{\mathrm{P}}$. This observation was confirmed for a wide range of models and parameter values through extensive numerical studies. It was also seen to hold for the endpoint-corrected estimators $A_{n,c}^{\mathrm{P}}$ and $A_{n,c}^{\mathrm{CFG}}$ originally proposed by Deheuvels (1991) and by Capéraà, Fougères and Genest (1997), respectively. It may be *conjectured*, therefore, that

$$\text{P–c} \xrightarrow{*} \text{CFG–c} \quad \text{and} \quad \text{P–opt} \xrightarrow{*} \text{CFG–opt}.$$

**5. Simulations.** A vast Monte Carlo study was used to confirm that the conclusions of Section 4 remain valid in finite-sample settings. For brevity, the results of a single experiment are reported here for illustration purposes.

Specifically, 5000 random samples of size $n = 100$ were generated from the $t$-EV copula with one degree of freedom and various values of $\rho$ chosen in such a way that the tail coefficient $2\{1 - A(0.5)\}$ ranges over the set $\{i/10 : i = 0, \dots, 10\}$. For each sample, the Hall–Tajvidi and the CFG estimators were computed when the margins are known and unknown. For each estimator, the empirical version of the mean integrated squared error,

$$\mathrm{MISE} = \mathrm{E}\left[\int_0^1 \{\hat{A}_n(t) - A(t)\}^2 \, dt\right],$$

was computed by averaging out over the 5000 samples.

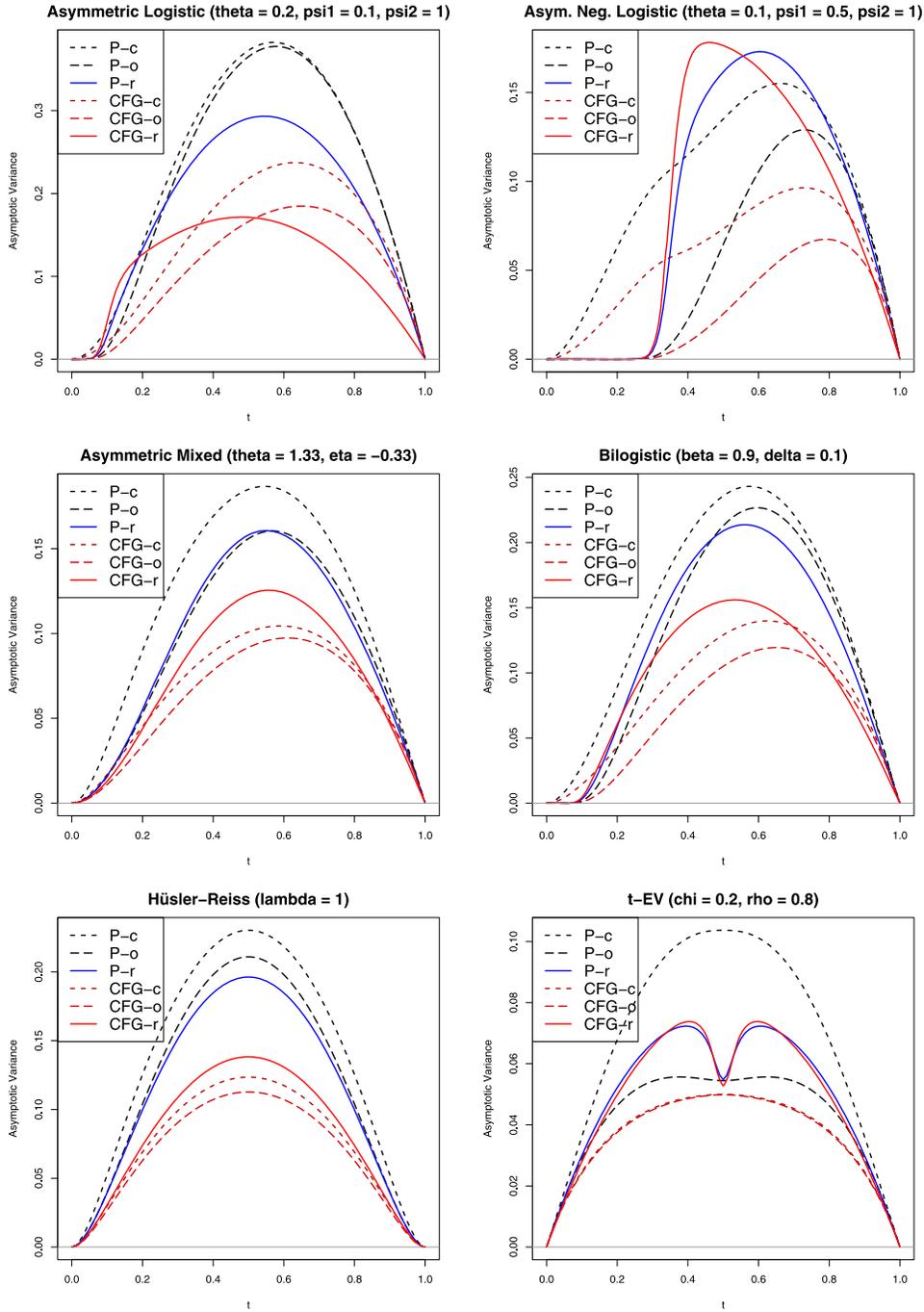

FIG. 3.  *Graph, as a function of $t$, of the asymptotic variances of the estimators $A_{n,\mathrm{c}}^{\mathrm{P}}(t)$,*
*$A_{n,\mathrm{opt}}^{\mathrm{P}}(t)$, $A_{n,\mathrm{r}}^{\mathrm{P}}(t)$, $A_{n,\mathrm{c}}^{\mathrm{CFG}}(t)$, $A_{n,\mathrm{opt}}^{\mathrm{CFG}}(t)$ and $A_{n,\mathrm{r}}^{\mathrm{CFG}}(t)$ for six extreme-value copula models.*

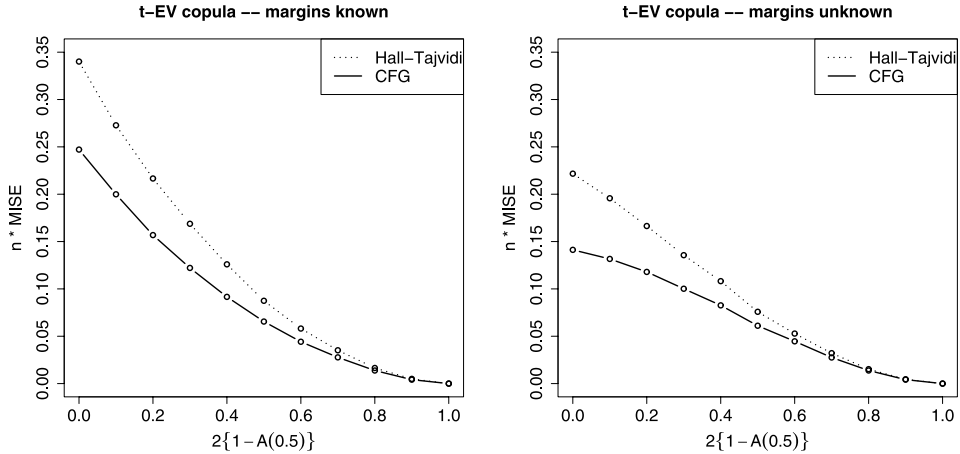

Fig. 4. MISE($\times 100$) of the estimators $A_n^{\mathrm{HT}}$ and $A_n^{\mathrm{CFG}}$ with correction $p(t) = 1 - t$, based on 5000 samples of size $n = 100$ from the $t$-EV copula with $\chi = 1$ degree of freedom and $\rho \in [-1, 1]$ chosen in such a way that $2\{1 - A(0.5)\} \in \{i/10 : i = 0, \ldots, 10\}$, when the margins are either known (left) or unknown (right).

In Figure 4, the normalized MISE (i.e., multiplied by $n$) is plotted as a function of $2\{1 - A(0.5)\}$. Similar results were obtained for many other extreme-value dependence models. In all cases, the CFG estimator was superior to the Hall–Tajvidi estimator. In addition, the use of ranks led to higher accuracy, although the gain diminished as the level of dependence increased. In the case of perfect positive dependence, $C(u, v) = u \wedge v$, the MISE of all estimators vanishes; indeed, it can be checked readily that their rate of convergence is then $o_p(n^{-1/2})$.

## 6. Conclusion.
This paper has proposed rank-based versions of several nonparametric estimators for the Pickands dependence function of a bivariate extreme-value copula. The new estimators were shown to be asymptotically normal and unbiased. Explicit formulas and consistent estimates for their asymptotic variances were also provided.

In contrast with the existing estimators, the new ones can be used regardless of whether the margins are known or unknown. It is interesting to note that the rank-based versions are generally more efficient than their uncorrected counterparts, even when the margins are known. In practice, however, only the endpoint corrected versions would be used. While the rank-based estimators no longer have a distinct advantage, they clearly remain competitive. The extensive numerical work presented herein also suggests that the CFG estimator is preferable to the Pickands estimator when (possibly optimal) endpoint corrections are applied to both of them.

A solution to this conjecture would clearly be of practical importance. In future work, it would also be of interest to assess more fully the impact of sample size, level and form of dependence (e.g., symmetry or lack thereof) on the precision of the various rank-based estimators proposed herein. Multivariate extensions could also be developed along similar lines as in Zhang, Wells and Peng (2008). Finally, the problem of constructing a truly semiparametric efficient estimator of a Pickands dependence function remains open.

## APPENDIX A: PROOF OF LEMMA 3.1

To express the Pickands estimator in the form (3.2), use Definition (2.2) and the substitution $u = e^{-x}$ to get

$$1/A_{n,\mathrm{r}}^{\mathrm{P}}(t) = \frac{1}{n}\sum_{i=1}^{n}\int_{0}^{\infty}\mathbf{1}\{\hat{\xi}_i(t) \geq x\}\,dx$$

$$= \frac{1}{n}\sum_{i=1}^{n}\int_{0}^{\infty}\mathbf{1}\{\hat{S}_i \geq (1-t)x, \hat{T}_i \geq tx\}\,dx$$

$$= \frac{1}{n}\sum_{i=1}^{n}\int_{0}^{\infty}\mathbf{1}\{\hat{U}_i \leq e^{-(1-t)x}, \hat{V}_i \leq e^{-tx}\}\,dx$$

$$= \frac{1}{n}\sum_{i=1}^{n}\int_{0}^{1}\mathbf{1}(\hat{U}_i \leq u^{1-t}, \hat{V}_i \leq u^t)\,\frac{du}{u} = \int_{0}^{1}\hat{C}_n(u^{1-t}, u^t)\,\frac{du}{u}.$$

In the case of the CFG estimator, the substitution $u = e^{-x}$ yields

$$\log(z) = \int_{1}^{\infty}\mathbf{1}(x \leq z)\,\frac{dx}{x} - \int_{0}^{1}\mathbf{1}(x > z)\,\frac{dx}{x}$$

$$= \int_{0}^{\infty}\{\mathbf{1}(x \leq z) - \mathbf{1}(x \leq 1)\}\,\frac{dx}{x}$$

$$= -\int_{0}^{1}[\mathbf{1}\{-\log(u) \leq z\} - \mathbf{1}(u \geq e^{-1})]\,\frac{du}{u\log u}.$$

Now, in light of Definition (2.3), one can write

$$-\gamma - \log A_{n,\mathrm{r}}^{\mathrm{CFG}}(t) = \frac{1}{n}\sum_{i=1}^{n}\log\hat{\xi}_i(t)$$

$$= \frac{1}{n}\sum_{i=1}^{n}\int_{0}^{1}[\mathbf{1}(u \geq e^{-1}) - \mathbf{1}\{-\log(u) \leq \hat{\xi}_i(t)\}]\,\frac{du}{u\log u}.$$

Further, note that $\hat{\xi}_i(t) \geq -\log(u)$ if and only if $-\log \hat{U}_i \geq -(1-t)\log(u)$ and $-\log \hat{V}_i \geq -t \log(u)$. The rest of the argument is along the same lines as in the case of the Pickands estimator.

## APPENDIX B: PROOF OF THEOREM 3.2

Before proceeding with the proof, recall that in Deheuvels (1979) the empirical copula is defined for all $u, v \in [0, 1]$ by

$$\hat{C}_n^{\mathrm{D}}(u, v) = H_n\{F_n^{\leftarrow}(u), G_n^{\leftarrow}(v)\}$$

in terms of the empirical marginal and joint distribution functions $F_n$, $G_n$ and $H_n$ of the sample $(X_1, Y_1), \ldots, (X_n, Y_n)$. Elementary calculations imply that, in the absence of ties (which occur with probability zero when the margins are continuous), the difference between $\hat{C}_n^{\mathrm{D}}$ and the empirical copula $\hat{C}_n$ in (3.1) is asymptotically negligible, namely,

$$(\text{B.1}) \qquad \sup_{(u,v)\in[0,1]^2} |\hat{C}_n^{\mathrm{D}}(u,v) - \hat{C}_n(u,v)| \leq \frac{4}{n}.$$

Let $C_n$ be the empirical distribution function of the (unobservable) random sample $(U_1, V_1), \ldots, (U_n, V_n)$, and set $\alpha_n = n^{1/2}(C_n - C)$. Further, define the (random) remainder term $R_n(u, v)$ implicitly by

$$(\text{B.2}) \quad \begin{aligned} &n^{1/2}\{\hat{C}_n(u,v) - C(u,v)\} \\ &= \alpha_n(u,v) - \dot{C}_1(u,v)\alpha_n(u,1) - \dot{C}_2(u,v)\alpha_n(1,v) + R_n(u,v). \end{aligned}$$

Now if $\hat{C}_n$ is replaced by $\hat{C}_n^{\mathrm{D}}$ in the left-hand side of (B.2), results in Stute [(1984), page 371] and Tsukahara [(2005), middle of page 359] imply that

$$(\text{B.3}) \qquad \sup_{(u,v)\in[0,1]} |R_n(u,v)| = O\{n^{-1/4}(\log n)^{1/2}(\log\log n)^{1/4}\}$$

almost surely as $n \to \infty$, provided that the second-order partial derivatives of $C$ exist and are continuous. This condition is automatically satisfied if $A$ is twice continuously differentiable, and hence in view of (B.1), the bound (B.3) remains valid for $R_n$ defined via $\hat{C}_n$.

PROOF OF THEOREM 3.2.

PICKANDS ESTIMATOR. Recall the empirical copula process $\mathbb{C}_n$ defined in (3.4), and for every $t \in [0, 1]$, let

$$\mathbb{B}_{n,\mathrm{r}}^{\mathrm{P}}(t) = n^{1/2}\{1/A_{n,\mathrm{r}}^{\mathrm{P}}(t) - 1/A(t)\}.$$

It will be shown that $\mathbb{B}_{n,r}^P \rightsquigarrow \mathbb{B} = -\mathbb{A}_r^P/A^2$ as $n \to \infty$, which implies that

$$\mathbb{A}_{n,r}^P = \frac{-A^2 \mathbb{B}_{n,r}^P}{1 + n^{1/2} \mathbb{B}_{n,r}^P} \rightsquigarrow \mathbb{A}_r^P$$

as a consequence of the functional version of Slutsky's lemma given in van der Vaart and Wellner [(1996), Examples 1.4.7, page 32].

Put $k_n = 2\log(n+1)$ and use identity (3.2) to see that

$$\mathbb{B}_{n,r}^P(t) = \int_0^1 \mathbb{C}_n(u^{1-t}, u^t) \frac{du}{u} = \int_0^\infty \mathbb{C}_n(e^{-s(1-t)}, e^{-st}) \, ds.$$

One can then write $\mathbb{B}_{n,r}^P(t) = I_{1,n}(t) + I_{2,n}(t)$, where for each $t \in [0,1]$,

$$I_{1,n}(t) = \int_{k_n}^\infty \mathbb{C}_n(e^{-s(1-t)}, e^{-st}) \, ds, \qquad I_{2,n}(t) = \int_0^{k_n} \mathbb{C}_n(e^{-s(1-t)}, e^{-st}) \, ds.$$

Note that the contribution of $I_{1,n}(t)$ is asymptotically negligible. For, if $s > k_n$, then $e^{-s(1-t)} \wedge e^{-st} \leq 1/(n+1)$, and hence

$$\mathbb{C}_n(e^{-s(1-t)}, e^{-st}) = -n^{1/2} C(e^{-s(1-t)}, e^{-st}) = -n^{1/2} e^{-sA(t)}.$$

Thus, for all $t \in [0,1]$,

$$(B.4) \quad |I_{1,n}(t)| = n^{1/2} \int_{k_n}^\infty e^{-sA(t)} \, ds \leq n^{1/2} \int_{k_n}^\infty e^{-s/2} \, ds = \frac{n^{1/2}}{n+1} \leq \frac{1}{n^{1/2}}.$$

Consequently, the asymptotic behavior of $\mathbb{B}_{n,r}^P$ is determined entirely by $I_{2,n}$. In turn, one can use Stute's representation (B.2) to write $I_{2,n} = J_{1,n} + \cdots + J_{4,n}$, where for each $t \in [0,1]$,

$$J_{1,n}(t) = \int_0^{k_n} \alpha_n(e^{-s(1-t)}, e^{-st}) \, ds,$$

$$J_{2,n}(t) = -\int_0^{k_n} \alpha_n(e^{-s(1-t)}, 1) \dot{C}_1(e^{-s(1-t)}, e^{-st}) \, ds,$$

$$J_{3,n}(t) = -\int_0^{k_n} \alpha_n(1, e^{-st}) \dot{C}_2(e^{-s(1-t)}, e^{-st}) \, ds,$$

$$J_{4,n}(t) = \int_0^{k_n} R_n(e^{-s(1-t)}, e^{-st}) \, ds.$$

From (B.3), the contribution of $J_{4,n}(t)$ becomes negligible as $n \to \infty$, because with probability one,

$$\sup_{t \in [0,1]} |J_{4,n}(t)| = O\{n^{-1/4}(\log n)^{3/2}(\log \log n)^{1/4}\}.$$

Accordingly, the identity $I_{2,n}(t) = J_{1,n}(t) + J_{2,n}(t) + J_{3,n}(t) + o(1)$ holds almost surely and uniformly in $t \in [0,1]$ when $n \to \infty$.

To complete the proof, it remains to show that the terms $J_{i,n}(t)$ with $i = 1, 2, 3$ have a suitable joint limit. To this end, fix $\omega \in (0, 1/2)$ and write $q_\omega(t) = t^\omega (1-t)^\omega$ for all $t \in [0, 1]$. Let also $\mathbb{G}_{n,\omega}(u, v) = \alpha_n(u, v)/q_\omega(u \wedge v)$ for all $u, v \in [0, 1]$. A simple substitution then shows that each $J_{i,n}(t)$ is a continuous functional of the weighted empirical process $\mathbb{G}_{n,\omega}$, namely,

$$(\text{B.5}) \qquad J_{1,n}(t) = \int_0^{k_n} \mathbb{G}_{n,\omega}(e^{-s(1-t)}, e^{-st}) K_1(s,t) \, ds,$$

$$(\text{B.6}) \qquad J_{2,n}(t) = -\int_0^{k_n} \mathbb{G}_{n,\omega}(e^{-s(1-t)}, 1) K_2(s,t) \, ds,$$

$$(\text{B.7}) \qquad J_{3,n}(t) = -\int_0^{k_n} \mathbb{G}_{n,\omega}(1, e^{-st}) K_3(s,t) \, ds$$

with $K_1(s,t) = q_\omega(e^{-s(1-t)} \wedge e^{-st})$, $K_2(s,t) = q_\omega(e^{-s(1-t)})\dot{C}_1(e^{-s(1-t)}, e^{-st})$ and $K_3(s,t) = q_\omega(e^{-st})\dot{C}_2(e^{-s(1-t)}, e^{-st})$ for all $s \in (0, \infty)$ and $t \in [0, 1]$.

The conclusion will then follow from Theorem G.1 and the continuous mapping theorem, provided that for $i = 1, 2, 3$, there exists an integrable function $K_i^* : (0, \infty) \to \mathbb{R}$ such that $K_i(s,t) \leq K_i^*(s)$ for all $s$ and $t$. For $K_1$, this is immediate because $K_1(s) \leq q_\omega(e^{-s/2}) \leq e^{-\omega s/2}$. For $K_2$, note that

$$\dot{C}_1(e^{-s(1-t)}, e^{-st}) = e^{-s\{A(t) - (1-t)\}} \mu(t),$$

where $\mu(t) = A(t) - t\dot{A}(t)$. Recalling also that $A(t) \geq t \vee (1-t)$, one finds

$$K_2(s,t) \leq \mu(t) e^{-s\{A(t) - (1-\omega)(1-t)\}} \leq \mu(t) e^{-\omega s/2}.$$

As a similar argument works for $K_3$, the proof is complete.

CFG ESTIMATOR. The argument mimics the proof pertaining to the Pickands estimator. To emphasize the parallel, the same notation is used and the presentation focuses on the changes.

In view of Lemma 3.1 and the functional version of Slutsky's lemma, the process to be studied is given for all $t \in [0, 1]$ by

$$\mathbb{B}_{n,\text{r}}^{\text{CFG}}(t) = n^{1/2}\{\log A_{n,\text{r}}^{\text{CFG}}(t) - \log A(t)\}$$

$$= \int_0^1 \mathbb{C}_n(u^{1-t}, u^t) \frac{du}{u \log u} = -\int_0^\infty \mathbb{C}_n(e^{-s(1-t)}, e^{-st}) \frac{ds}{s}.$$

This process can be decomposed as $-(I_{1,n} + I_{2,n} + I_{3,n})$, where

$$I_{1,n}(t) = \int_{k_n}^\infty \mathbb{C}_n(e^{-s(1-t)}, e^{-st}) \frac{ds}{s},$$

$$I_{2,n}(t) = \int_{\ell_n}^{k_n} \mathbb{C}_n(e^{-s(1-t)}, e^{-st}) \frac{ds}{s},$$

$$I_{3,n}(t) = \int_0^{\ell_n} \mathbb{C}_n(e^{-s(1-t)}, e^{-st}) \frac{ds}{s}$$

with $k_n = 2\log(n+1)$ as above and $\ell_n = 1/(n+1)$.

Arguing as in (B.4), one sees that $|I_{1,n}| \leq n^{-1/2}$. Similarly, $I_{3,n}$ is negligible asymptotically. For, if $s \in (0, \ell_n)$ and $t \in [0,1]$, one has $e^{-s(1-t)} \wedge e^{-st} \geq e^{-1/(n+1)} > 1/(n+1)$, and then $\hat{C}_n(e^{-s(1-t)}, e^{-st}) = 1$, so that

$$|\mathbb{C}_n(e^{-s(1-t)}, e^{-st})| \leq n^{1/2}\{1 - e^{-sA(t)}\} \leq n^{1/2} s A(t).$$

Therefore, $|I_{3,n}| \leq n^{1/2}\ell_n \leq n^{-1/2}$.

As a result, the asymptotic behavior of $\mathbb{B}_{n,\mathrm{r}}^{\mathrm{CFG}}$ is determined entirely by $I_{2,n}$. Using Stute's representation (B.2), one may then write $I_{2,n} = J_{1,n} + \cdots + J_{4,n}$, where for all $t \in [0,1]$,

$$J_{1,n}(t) = \int_{\ell_n}^{k_n} \alpha_n(e^{-s(1-t)}, e^{-st}) \frac{ds}{s},$$

$$J_{2,n}(t) = -\int_{\ell_n}^{k_n} \alpha_n(e^{-s(1-t)}, 1)\dot{C}_1(e^{-s(1-t)}, e^{-st}) \frac{ds}{s},$$

$$J_{3,n}(t) = -\int_{\ell_n}^{k_n} \alpha_n(1, e^{-st})\dot{C}_2(e^{-s(1-t)}, e^{-st}) \frac{ds}{s},$$

$$J_{4,n}(t) = \int_{\ell_n}^{k_n} R_n(e^{-s(1-t)}, e^{-st}) \frac{ds}{s}.$$

Again, the term $J_{4,n}$ is negligible asymptotically because as $n \to \infty$,

$$|J_{4,n}(t)| \leq \log(k_n/\ell_n) \sup_{(u,v)\in[0,1]^2} |R_n(u,v)|$$

$$= O\{n^{-1/4}(\log n)^{3/2}(\log\log n)^{1/4}\}$$

almost surely and uniformly in $t \in [0,1]$. As for $J_{1,n}$, $J_{2,n}$ and $J_{3,n}$, they admit the same representations as (B.5)–(B.7), except that in each case, the integration is limited to the interval $(\ell_n, k_n)$. For $s \in [1, \infty)$, the same upper bounds $K_1^*$, $K_2^*$, $K_3^*$ apply, and they have already been shown to be integrable on this domain. As for the integrability on $(0,1)$, it follows from the additional bound $|1 - e^{-s(1-t)} \wedge e^{-st}|^\omega \leq s^\omega$.   $\square$

## APPENDIX C: PROOF OF PROPOSITION 3.3

The proofs rely on the fact that for all $u, v, t \in [0,1]$,

$$\sigma(u,v;t) = \sigma_0(u,v;t) + (uv)^{A(t)}\left\{\sum_{\ell=1}^{4} \sigma_\ell(u,v;t) - \sum_{\ell=5}^{8} \sigma_\ell(u,v;t)\right\},$$

where $\sigma_0$ is given by (4.1) and

$$\sigma_1(u,v;t) = (u^{t-1} \wedge v^{t-1} - 1)\mu^2(t),$$

$$\sigma_2(u,v;t) = (u^{-t} \wedge v^{-t} - 1)\nu^2(t),$$

$$\sigma_3(u,v;t) = \{u^{t-1}v^{-t}C(u^{1-t},v^t) - 1\}\mu(t)\nu(t),$$

$$\sigma_4(u,v;t) = \{u^{-t}v^{t-1}C(v^{1-t},u^t) - 1\}\mu(t)\nu(t),$$

$$\sigma_5(u,v;t) = \{u^{-A(t)}v^{t-1}C(u^{1-t} \wedge v^{1-t}, u^t) - 1\}\mu(t),$$

$$\sigma_6(u,v;t) = \{u^{t-1}v^{-A(t)}C(u^{1-t} \wedge v^{1-t}, v^t) - 1\}\mu(t),$$

$$\sigma_7(u,v;t) = \{u^{-A(t)}v^{-t}C(u^{1-t}, u^t \wedge v^t) - 1\}\nu(t),$$

$$\sigma_8(u,v;t) = \{u^{-t}v^{-A(t)}C(v^{1-t}, u^t \wedge v^t) - 1\}\nu(t).$$

As $A_{n,\mathrm{r}}^{\mathrm{P}}(0) = A_{n,\mathrm{r}}^{\mathrm{P}}(1) = 1$ and $A_{n,\mathrm{r}}^{\mathrm{CFG}}(0) = A_{n,\mathrm{r}}^{\mathrm{CFG}}(1) = 1$ by construction, it may be assumed without loss of generality that $t \in (0,1)$. It is immediate from the work of Segers (2007) that

$$\int_0^1 \int_0^1 \sigma_0(u,v;t)\, \frac{du}{u}\, \frac{dv}{v} = \frac{1}{A^2(t)}$$

and

$$\int_0^1 \int_0^1 \sigma_0(u,v;t)\, \frac{du}{u\log(u)}\, \frac{dv}{v\log(v)} = L_2(1)$$

for the Pickands and the CFG estimators, respectively. Various symmetries help to reduce the computation of the remaining terms from eight to three for both estimators. In particular, note that $\sigma_i(u,v;t) = \sigma_{i+1}(v,u;t)$ for $i = 3,5,7$ and all $t \in [0,1]$. Furthermore,

$$\sigma_2(u,v;t) = \overline{\sigma_1(u,v;t)}, \qquad \sigma_7(u,v;t) = \overline{\sigma_5(u,v;t)},$$

where a bar over a function means that all instances of $A$, $t$, $\mu$ and $\nu$ in it should be replaced by $1-A$, $1-t$, $1-\mu$ and $1-\nu$, respectively. Thus, if

$$f(u,v) = \begin{cases} f_P(u,v) = (uv)^{A(t)-1}, \\ f_{\mathrm{CFG}}(u,v) = (uv)^{A(t)-1}/\{\log(u)\log(v)\}, \end{cases}$$

one has $f(v,u) = f(u,v) = \overline{f(u,v)}$, as well as

$$\sum_{i=1}^4 \int_0^1 \int_0^1 \sigma_i(u,v;t) f(u,v)\, du\, dv$$

$$= \int_0^1 \int_0^1 \{\sigma_1(u,v;t) + \overline{\sigma_1(u,v;t)} + 2\sigma_3(u,v;t)\} f(u,v)\, du\, dv$$

and

$$\sum_{i=5}^{8} \int_0^1 \int_0^1 \sigma_i(u,v;t) f(u,v)\, du\, dv$$

$$= \int_0^1 \int_0^1 \{2\sigma_5(u,v;t) + 2\overline{\sigma_5(u,v;t)}\} f(u,v)\, du\, dv.$$

Each of the relevant parts is computed in turn.

INTEGRALS INVOLVING $\sigma_1$. For the Pickands estimator,

$$\int_0^1 \int_0^1 (u \wedge v)^{t-1}(uv)^{A(t)-1}\, du\, dv$$

(C.1)
$$= 2\int_0^1 \int_0^v u^{A(t)-1} v^{A(t)+t-2}\, du\, dv$$

$$= \frac{2}{A(t)} \int_0^1 v^{2A(t)+t-2}\, dv = \frac{2}{A(t)\{2A(t)+t-1\}}.$$

Consequently,

$$\int_0^1 \int_0^1 \sigma_1(u,v;t) f_P(u,v)\, du\, dv = \frac{\mu^2(t)}{A(t)} \left\{ \frac{2}{2A(t)+t-1} - \frac{1}{A(t)} \right\}.$$

For the CFG estimator,

$$\int_0^1 \int_0^1 \{(u \wedge v)^{t-1} - 1\}(uv)^{A(t)} \frac{du}{u\log(u)} \frac{dv}{v\log(v)}$$

(C.2)
$$= 2\int_0^1 \int_u^1 (v^{t-1} - 1)(uv)^{A(t)} \frac{dv}{v\log(v)} \frac{du}{u\log(u)}.$$

Use the substitution $v = u^x$ to rewrite this expression as

$$-2\int_0^1 \int_0^1 (u^{(t-1)x} - 1)u^{(1+x)A(t)-1} \frac{du}{\log(u)} \frac{dx}{x}$$

$$= -2\int_0^1 \log\left\{ \frac{(1+x)A(t)+(t-1)x}{(1+x)A(t)} \right\} \frac{dx}{x}$$

$$= -2\int_0^1 \log[1 + \{1 - 1/A_2(t)\}x] \frac{dx}{x} + 2\int_0^1 \log(1+x) \frac{dx}{x}$$

$$= 2L_2\{-1 + 1/A_2(t)\} + L_2(1).$$

Therefore,

$$\int_0^1 \int_0^1 \sigma_1(u,v;t) f_{\mathrm{CFG}}(u,v)\, du\, dv = \mu^2(t)[2L_2\{-1 + 1/A_2(t)\} + L_2(1)].$$

INTEGRALS INVOLVING $\sigma_3$. First, consider the Pickands estimator. The substitutions $u^{1-t} = x$ and $v^t = y$ yield

$$\int_0^1 \int_0^1 u^{t-1} v^{-t} C(u^{1-t}, v^t)(uv)^{A(t)-1}\, du\, dv$$

$$= \frac{1}{t(1-t)} \int_0^1 \int_0^1 \frac{C(x,y)}{xy}(x^{1/(1-t)} y^{1/t})^{A(t)-1} x^{1/(1-t)-1} y^{1/t-1}\, dx\, dy$$

$$= \frac{1}{t(1-t)} \int_0^1 \int_0^1 C(x,y) x^{A_2(t)-2} y^{A_1(t)-2}\, dx\, dy.$$

Next, use the substitutions $x = w^{1-s}$ and $y = w^s$. Note that $w = xy \in (0,1]$, $s = \log(y)/\log(xy) \in [0,1]$, $C(x,y) = w^{A(s)}$ and the Jacobian of the transformation is $-\log(w)$. The above integral then becomes

$$-\frac{1}{t(1-t)} \int_0^1 \int_0^1 w^{A(s)+(1-s)\{A_2(t)-2\}+s\{A_1(t)-2\}} \log(w)\, dw\, ds$$

$$= \frac{1}{t(1-t)} \int_0^1 [A(s) + (1-s)\{A_2(t)-2\} + s\{A_1(t)-2\} + 1]^{-2}\, ds.$$

With $A(s)$ in the integrand, no further simplification is possible. Thus,

$$\int_0^1 \int_0^1 \sigma_3(u,v;t) f_P(u,v)\, du\, dv$$

$$= \frac{\mu(t)\nu(t)}{A^2(t)} - \frac{\mu(t)\nu(t)}{t(1-t)} \int_0^1 \{A(s) + sA_1(t) + (1-s)A_2(t) - 1\}^{-2}\, ds.$$

The same substitutions are used for the CFG estimator. They yield

$$\int_0^1 \int_0^1 \left\{ \frac{C(u^{1-t}, v^t)}{u^{1-t} v^t} - 1 \right\} (uv)^{A(t)} \frac{du}{u \log(u)} \frac{dv}{v \log(v)}$$

$$= \int_0^1 \int_0^1 \left\{ \frac{C(x,y)}{xy} - 1 \right\} x^{A_2(t)} y^{A_1(t)} \frac{dx}{x \log(x)} \frac{dy}{y \log(y)}$$

$$= -\int_0^1 \int_0^1 (w^{-1+A(s)} - 1) w^{\{(1-s)/(1-t)+s/t\}A(t)-1} \frac{dw}{\log(w)} \frac{ds}{s(1-s)}$$

$$= -\int_0^1 \log\left[ \frac{\{(1-s)/(1-t)+s/t\}A(t) - 1 + A(s)}{\{(1-s)/(1-t)+s/t\}A(t)} \right] \frac{ds}{s(1-s)}.$$

Hence,

$$\int_0^1 \int_0^1 \sigma_3(u,v;t) f_{\text{CFG}}(u,v)\, du\, dv$$

$$= -\mu(t)\nu(t) \int_0^1 \log\left\{ 1 - \frac{t(1-t)}{A(t)} \frac{1-A(s)}{t(1-s)+(1-t)s} \right\} \frac{ds}{s(1-s)}.$$

INTEGRALS INVOLVING $\sigma_5$. For the Pickands estimator, one gets

$$\int_0^1 \int_0^1 u^{-A(t)} v^{t-1} C\{(u \wedge v)^{1-t}, u^t\}(uv)^{A(t)} \frac{du}{u} \frac{dv}{v}$$

$$= \int_0^1 \int_0^v v^{t-1}(uv)^{A(t)-1} \, du \, dv$$

$$+ \int_0^1 \int_v^1 u^{-A(t)} v^{t-1} C(v^{1-t}, u^t)(uv)^{A(t)-1} \, du \, dv.$$

The first integral on the right was computed in (C.1). For the second, collect powers in $u$ and $v$ and use the substitutions $v^{1-t} = x$ and $u^t = y$. This yields

$$\frac{1}{t(1-t)} \int_0^1 \int_0^{y^{-1+1/t}} x^{A_2(t)-2} y^{-1} C(x, y) \, dx \, dy.$$

Then, make the same substitutions $x = w^{1-s}$, $y = w^s$ that were used for $\sigma_3$. As the constraint $x < y^{-1+1/t}$ reduces to $s < t$, the integral becomes

$$-\frac{1}{t(1-t)} \int_0^t \int_0^1 w^{(1-s)\{A_2(t)-2\}+A(s)-s} \log(w) \, dw \, ds$$

$$= \frac{1}{t(1-t)} \int_0^t [(1-s)\{A_2(t)-2\} + A(s) - s + 1]^{-2} \, ds.$$

Thus,

$$\int_0^1 \int_0^1 \sigma_5(u, v; t) f_P(u, v) \, du \, dv$$

$$= \frac{\mu(t)}{A(t)\{2A(t)+t-1\}} - \frac{\mu(t)}{A^2(t)}$$

$$+ \frac{\mu(t)}{t(1-t)} \int_0^t [A(s) + (1-s)\{A_2(t)-1\}]^{-2} \, ds.$$

Finally, turning to the CFG estimator, one must evaluate

$$\int_0^1 \int_0^v (v^{t-1} - 1)(uv)^{A(t)} \frac{du}{u \log(u)} \frac{dv}{v \log(v)}$$

$$+ \int_0^1 \int_v^1 \{u^{-A(t)} v^{t-1} C(v^{1-t}, u^t) - 1\}(uv)^{A(t)} \frac{du}{u \log(u)} \frac{dv}{v \log(v)}.$$

The first integral on the right was computed in (C.2). For the second, the same substitutions are used as for the Pickands estimator. This yields

$$\int_0^1 \int_0^{y^{-1+1/t}} \{x^{-1} y^{-A_1(t)} C(x, y) - 1\} x^{A_2(t)} y^{A_1(t)} \frac{dx}{x \log(x)} \frac{dy}{y \log y}$$

$$= -\int_0^t \int_0^1 (w^{s-1-sA_1(t)+A(s)} - 1)w^{\{(1-s)/(1-t)+s/t\}A(t)-1}$$

$$\times \frac{dw}{\log(w)} \frac{ds}{s(1-s)}$$

$$= -\int_0^t \log\left[\frac{\{(1-s)/(1-t)+s/t\}A(t)+s-1-s/tA(t)+A(s)}{\{(1-s)/(1-t)+s/t\}A(t)}\right]$$

$$\times \frac{ds}{s(1-s)}.$$

Therefore,

$$\int_0^1 \int_0^1 \sigma_5(u,v;t) f_{\mathrm{CFG}}(u,v)\, du\, dv$$

$$= \mu(t)L_2\{-1+1/A_2(t)\} + \mu(t)L_2(1)/2$$

$$- \mu(t)\int_0^t \log\left[1 - \frac{t(1-t)}{A(t)} \frac{1-A(s)+s\{A_1(t)-1\}}{t(1-s)+(1-t)s}\right] \frac{ds}{s(1-s)}.$$

It then suffices to assemble the various terms to conclude.

## APPENDIX D: PROOF OF COROLLARY 3.4

When $C(u,v) = uv$ for all $u,v \in [0,1]$, one has $A(t) = \mu(t) = \nu(t) = 1$ for all $t \in [0,1]$. Upon substitution, one gets

$$\sigma(u,v;t) = (u^{1-t} \wedge v^{1-t} - u^{1-t}v^{1-t})(u^t \wedge v^t - u^t v^t),$$

which simplifies to $\sigma(u,v;t) = u(1-v^{1-t})(1-v^t)$ for arbitrary $u,v \in (0,1)$ with $u < v$. Thus, by symmetry,

$$\mathrm{var}\, \mathbb{A}_{\mathrm{r}}^{\mathrm{P}}(t) = 2\int_0^1 (1-v^{1-t})(1-v^t)\, dv = \frac{3t(1-t)}{(2-t)(1+t)}$$

for all $t \in [0,1]$. For the CFG estimator, the substitution $v = u^x$ yields

$$\mathrm{var}\, \mathbb{A}_{\mathrm{r}}^{\mathrm{CFG}}(t) = 2\int_0^1 \int_u^1 (1-v^{1-t})(1-v^t)\, \frac{dv}{v\log(v)} \frac{du}{\log(u)}$$

$$= -2\int_0^1 \int_0^1 (1-u^{(1-t)x})(1-u^{xt})\, \frac{du}{\log(u)} \frac{dx}{x}$$

$$= 2\int_0^1 \left[\log\{(1-t)x+1\} - \log\left\{\frac{(1-t)x+xt+1}{1+xt}\right\}\right] \frac{dx}{x}.$$

The latter decomposes into a sum of three integrals, namely,

$$-2\int_0^1 \log(1+x)\, \frac{dx}{x} + 2\int_0^{1-t} \log(1+x)\, \frac{dx}{x} + 2\int_0^t \log(1+x)\, \frac{dx}{x},$$

whence the conclusion.

## APPENDIX E: PROOF OF PROPOSITION 4.1

COMPARISON OF THE PICKANDS ESTIMATORS. In view of Theorem 3.1 of Segers (2007) and relation (4.2) with $a(t) = t$ and $b(t) = 1 - t$, one has

$$\operatorname{var} \mathbb{A}_c^{\mathrm{P}}(t) - \operatorname{var} \mathbb{A}_u^{\mathrm{P}}(t) = (1-t)^2 \operatorname{var} \eta(0) + t^2 \operatorname{var} \eta(1) - 2(1-t) \operatorname{cov}\{\eta(0), \eta(t)\}$$
$$- 2t \operatorname{cov}\{\eta(t), \eta(1)\} + 2t(1-t) \operatorname{cov}\{\eta(0), \eta(1)\},$$

where $\eta$ denotes a zero-mean Gaussian process on $[0, 1]$ whose covariance function is defined for all $0 \le s \le t \le 1$ by

$$\operatorname{cov}\{\eta(s), \eta(t)\} = \frac{s}{t} \frac{1}{A^2(s)} + \frac{1-t}{1-s} \frac{1}{A^2(t)} + \frac{1}{(1-s)t} \int_s^t \frac{dw}{A^2(w)} - \frac{1}{A(s)A(t)}.$$

Upon substitution and simplification, $\operatorname{var} \mathbb{A}_c^{\mathrm{P}}(t) - \operatorname{var} \mathbb{A}_u^{\mathrm{P}}(t)$ thus reduces to

$$(\mathrm{E}.1) \quad \begin{aligned} &-1 + \frac{2}{A(t)} + 2\{(1-t)^2 + t^2\}\left\{1 - \frac{1}{A^2(t)}\right\} \\ &- 2(1-t)\left(\frac{1}{t} - t\right) \int_0^t \frac{dw}{A^2(w)} - 2t\left\{\frac{1}{1-t} - (1-t)\right\} \int_t^1 \frac{dw}{A^2(w)}. \end{aligned}$$

Because $A$ is convex, however, one knows that for all $t \in [0, 1]$,

$$\int_0^t \frac{dw}{A^2(w)} \ge \frac{1}{A(t)} \quad \text{and} \quad \int_t^1 \frac{dw}{A^2(w)} \ge \frac{1-t}{A(t)}.$$

Using these inequalities and the fact that the third summand in (E.1) is negative for all $t \in [0, 1]$, one gets

$$\operatorname{var} \mathbb{A}_c^{\mathrm{P}}(t) - \operatorname{var} \mathbb{A}_u^{\mathrm{P}}(t) \le -1 + 2[1 - (1-t)(1-t^2) - t\{1 - (1-t)^2\}]\frac{1}{A(t)}$$

$$= -1 + 2t(1-t)\frac{1}{A(t)}$$

and this upper bound is negative because $A(t) \ge t \vee (1 - t)$ for all $t \in [0, 1]$. Thus, the argument is complete.

COMPARISON OF THE CFG ESTIMATORS. Theorem 4.2 of Segers (2007) and relation (4.2) with $a(t) = t$ and $b(t) = 1 - t$ imply that

$$\operatorname{var} \mathbb{A}_c^{\mathrm{CFG}}(t) - \operatorname{var} \mathbb{A}_u^{\mathrm{CFG}}(t)$$
$$= (1-t)^2 \operatorname{var} \zeta(0) + t^2 \operatorname{var} \zeta(1) - 2(1-t) \operatorname{cov}\{\zeta(0), \zeta(t)\}$$
$$- 2t \operatorname{cov}\{\zeta(t), \zeta(1)\} + 2t(1-t) \operatorname{cov}\{\zeta(0), \zeta(1)\},$$

where $\zeta$ denotes a zero-mean Gaussian process on $[0,1]$ whose covariance function is defined for all $0 \le s \le t \le 1$ by

$$\operatorname{cov}\{\zeta(s), \zeta(t)\} = -\int_0^s \log(w)\, \frac{dw}{1-w} - \log(t)\log(1-s) - \int_t^1 \log(1-w)\, \frac{dw}{w}$$

$$+ \log\left(\frac{t}{s}\right)\log A(s) + \log\left(\frac{1-s}{1-t}\right)\log A(t)$$

$$+ \frac{1}{2}\left\{\log\frac{A(s)}{A(t)}\right\}^2 - \int_s^t \frac{\log A(w)}{w(1-w)}\, dw.$$

Upon substitution and simplification, $\operatorname{var}\mathbb{A}_{\mathrm{c}}^{\mathrm{CFG}}(t) - \operatorname{var}\mathbb{A}_{\mathrm{u}}^{\mathrm{CFG}}(t)$ becomes

$$\{(1-t)^2 + t^2\}\frac{\pi^2}{6} + 2(1-t)\int_t^1 \log(1-w)\, \frac{dw}{w} + 2t\int_0^t \log(w)\, \frac{dw}{1-w}$$

$$- \{\log A(t)\}^2 + 2\{(1-t)\log(1-t) + t\log(t)\}\log A(t)$$

$$+ 2(1-t)\int_0^t \frac{\log A(w)}{w(1-w)}\, dw + 2t\int_t^1 \frac{\log A(w)}{w(1-w)}\, dw$$

$$- 2t(1-t)\int_0^1 \frac{\log A(w)}{w(1-w)}\, dw.$$

Omitting the term $-\{\log A(t)\}^2$ and using the elementary inequalities

$$\frac{\pi^2}{6}(1-t) \le -\int_t^1 \log(1-w)\, \frac{dw}{w}, \qquad \frac{\pi^2}{6}t \le -\int_0^t \log(w)\, \frac{dw}{1-w},$$

one can see that an upper bound on $\operatorname{var}\mathbb{A}_{\mathrm{c}}^{\mathrm{CFG}}(t) - \operatorname{var}\mathbb{A}_{\mathrm{u}}^{\mathrm{CFG}}(t)$ is given by

$$(1-t)\int_t^1 \log(1-w)\, \frac{dw}{w} + t\int_0^t \log(w)\, \frac{dw}{1-w}$$

$$+ 2\{(1-t)\log(1-t) + t\log(t)\}\log A(t)$$

$$+ 2(1-t)^2 \int_0^t \frac{\log A(w)}{w(1-w)}\, dw + 2t^2 \int_t^1 \frac{\log A(w)}{w(1-w)}\, dw.$$

Partial integration and the fact that $A(t) \ge 1-t$ for all $t \in [0,1]$ imply that

$$\int_0^t \log(w)\frac{dw}{1-w} = -\log(t)\log(1-t) + \int_0^t \log(1-w)\, \frac{dw}{w}$$

$$\le -\log(t)\log(1-t) \le -\log(t)\log A(t).$$

Similarly, the fact that $A(t) \ge t$ for all $t \in [0,1]$ yields

$$\int_t^1 \log(1-w)\, \frac{dw}{w} \le -\log(1-t)\log A(t).$$

Therefore, a more conservative upper bound on $\operatorname{var} \mathbb{A}_c^{\mathrm{CFG}}(t) - \operatorname{var} \mathbb{A}_u^{\mathrm{CFG}}(t)$ is given by

$$\{(1-t)\log(1-t) + t\log(t)\}\log A(t)$$
$$+ 2(1-t)^2 \int_0^t \frac{\log A(w)}{w(1-w)}\,dw + 2t^2 \int_t^1 \frac{\log A(w)}{w(1-w)}\,dw.$$

Now $A(t) \in [1/2, 1]$ and hence $2\{A(t) - 1\} \le \log A(t) \le A(t) - 1$ for all $t \in [0, 1]$. Consequently, an even more conservative bound is given by

$$2\{(1-t)\log(1-t) + t\log(t)\}\{A(t) - 1\}$$
$$+ 2(1-t)^2 \int_0^t \frac{A(w) - 1}{w(1-w)}\,dw + 2t^2 \int_t^1 \frac{A(w) - 1}{w(1-w)}\,dw.$$

Calling on the convexity of $A$, however, one can show that $\{A(w) - 1\}/w \le \{A(t) - 1\}/t$ for all $w \in [0, t]$ while $\{A(w) - 1\}/(1-w) \le \{A(t) - 1\}/(1-t)$ for any $w \in [t, 1]$. This leads to the final upper bound, namely,

$$2\frac{1 - 2t}{t(1-t)}\{t^2\log(t) - (1-t)^2\log(1-t)\}\{A(t) - 1\}.$$

As the latter is easily checked to be at most zero, the proof is complete.

## APPENDIX F: PROOF OF PROPOSITION 4.2

In view of the expressions for $\sigma(u, v; t)$ and $\sigma_0(u, v; t)$ given in Appendix C, it suffices to show that for all $u, v, t \in [0, 1]$,

$$\sum_{\ell=1}^4 \sigma_\ell(u, v; t) \le \sum_{\ell=5}^8 \sigma_\ell(u, v; t).$$

If $u \le v$, then

$$\sigma_1 \le \sigma_5, \qquad \sigma_2 \le \sigma_7, \qquad \sigma_3 \le \sigma_6, \qquad \sigma_4 \le \sigma_8$$

for all $t \in [0, 1]$, while if $v \le u$, then

$$\sigma_1 \le \sigma_6, \qquad \sigma_2 \le \sigma_8, \qquad \sigma_3 \le \sigma_7, \qquad \sigma_4 \le \sigma_5$$

for all $t \in [0, 1]$. Each of these inequalities is an easy consequence of the following inequalities, which are valid for every extreme-value copula $C$, associated dependence function $A$ and real numbers $u, v, t \in [0, 1]$:

$$uv \le C(u, v), \qquad t \vee (1 - t) \le A(t) \le 1, \qquad 0 \le \mu(t) \le 1, \qquad 0 \le \nu(t) \le 1.$$

The bounds on $\mu(t)$ and $\nu(t)$ stem from the fact that $A$ is convex, $A(0) = A(1) = 1$ and $\dot{A}(t) \in [-1, 1]$ for every $t \in (0, 1)$.

## APPENDIX G: WEIGHTED BIVARIATE EMPIRICAL PROCESSES

Let $(U_1, V_1), \ldots, (U_n, V_n)$ be a random sample from an arbitrary bivariate copula $C$ and for all $u, v \in [0, 1]$, define

$$C_n(u, v) = \frac{1}{n} \sum_{i=1}^{n} \mathbf{1}(U_i \le u, V_i \le v).$$

The purpose of this appendix is to characterize the asymptotic behavior of a weighted version of the empirical process $\alpha_n = n^{1/2}(C_n - C)$. Specifically, fix $\omega \ge 0$, and for every $x \in [0, 1]$, let $q_\omega(t) = t^\omega (1 - t)^\omega$ and

$$\mathbb{G}_{n,\omega}(u, v) = \begin{cases} \dfrac{\alpha_n(u, v)}{q_\omega(u \wedge v)}, & \text{if } u \wedge v \in (0, 1), \\ 0, & \text{if } u = 0 \text{ or } v = 0 \text{ or } (u, v) = (1, 1). \end{cases}$$

The following result, which may be of independent interest, gives the weak limit of the weighted process $\mathbb{G}_{n,\omega}$ in the space $\ell^\infty([0, 1]^2)$ of bounded, real-valued functions on $[0, 1]^2$ equipped with the topology of uniform convergence. Weak convergence is understood in the sense of Hoffman-Jørgensen [van der Vaart and Wellner (1996), Section 1.5].

THEOREM G.1. *For every* $\omega \in [0, 1/2)$, *the process* $\mathbb{G}_{n,\omega}$ *converges weakly in* $\ell^\infty([0, 1]^2)$ *to a centered Gaussian process* $\mathbb{G}_\omega$ *with continuous sample paths such that* $\mathbb{G}_\omega(u, v) = 0$ *if* $u = 0$, $v = 0$ *or* $(u, v) = (1, 1)$, *while*

$$\text{cov}\{\mathbb{G}_\omega(u, v), \mathbb{G}_\omega(u', v')\} = \frac{C(u \wedge u', v \wedge v') - C(u, v)C(u', v')}{q_\omega(u \wedge v)q_\omega(u' \wedge v')},$$

*if* $u \wedge v \in (0, 1)$ *and* $u' \wedge v' \in (0, 1)$.

The proof of this result relies on the theory of empirical processes detailed in van der Vaart and Wellner (1996), whose notation is adopted. Let

$$\mathbb{E} = \{(u, v) \in [0, 1]^2 : 0 < u \wedge v < 1\} = (0, 1]^2 \setminus \{(1, 1)\}.$$

For fixed $(u, v) \in \mathbb{E}$, define the mapping $f_{(u,v)} : \mathbb{E} \to \mathbb{R}$ by

$$(s, t) \mapsto f_{(u,v)}(s, t) = \frac{\mathbf{1}_{(0,u] \times (0,v]}(s, t) - C(u, v)}{q(u \wedge v)}$$

and consider the class $\mathcal{F} = \{f_{(u,v)} : (u, v) \in \mathbb{E}\} \cup \{0\}$, where $0$ denotes the function vanishing everywhere on $\mathbb{E}$.

Finally, let $P$ be the probability distribution on $\mathbb{E}$ corresponding to $C$ and denote by $\mathbb{P}_n$ the empirical measure of the sample $(U_1, V_1), \ldots, (U_n, V_n)$. The following lemma is instrumental in the proof of Theorem G.1. It pertains to the asymptotic behavior of the process $\mathbb{G}_n$, where for each $f \in \mathcal{F}$,

$$\mathbb{G}_n f = n^{1/2}(\mathbb{P}_n f - Pf)$$

with

$$\mathbb{P}_n f = \frac{1}{n} \sum_{i=1}^{n} f(U_i, V_i), \qquad Pf = \iint_{[0,1]^2} f(u,v)\, dC(u,v).$$

LEMMA G.2.   *The collection $\mathcal{F}$ is a $P$-Donsker class, that is, there exists a $P$-Brownian bridge $\mathbb{G}$ such that $\mathbb{G}_n \rightsquigarrow \mathbb{G}$ as $n \to \infty$ in $\ell^\infty(\mathcal{F})$.*

PROOF.   It is enough to check the conditions of Dudley and Koltchinskii (1994) reported in Theorem 2.6.14 of van der Vaart and Wellner (1996), namely:

(a)  $\mathcal{F}$ is a VC-major class.
(b)  There exists $F\colon \mathbb{E} \to \mathbb{R}$ such that $|f| \le F$ pointwise for all $f \in \mathcal{F}$ and

$$\int_0^\infty \{P(F > x)\}^{1/2}\, dx < \infty,$$

where $P(F > x) = P\{(s,t) \in \mathbb{E} : F(s,t) > x\}$.
(c)  $\mathcal{F}$ is pointwise separable.

To prove (a), one must check that the class of subsets of $\mathbb{E}$ given by $\{(u',v') \in \mathbb{E} : f(u',v') > t\}$ with $f$ ranging over $\mathcal{F}$ and $t$ over $\mathbb{R}$ forms a Vapnik–Červonenkis (VC) class of sets; see page 145 in van der Vaart and Wellner (1996). Noting that each $f$ in $\mathcal{F}$ can take only two values, one can see that this class of subsets coincides with the family of intervals $(0, u] \times (0, v]$ with $(u, v)$ ranging over $\mathbb{E}$. By Example 2.6.1 on page 135 of van der Vaart and Wellner (1996), the latter is indeed a VC-class.

To prove (b), define $F\colon \mathbb{E} \to \mathbb{R}$ at every $(s, t) \in \mathbb{E}$ by

$$F(s,t) = 2\{s^{-\omega} \vee t^{-\omega} \vee (1-s)^{-\omega} \vee (1-t)^{-\omega}\}.$$

Given that the marginal distributions of $C$ are uniform on $[0, 1]$, one has

(G.1)   $C(u,v) \le u \wedge v$   and   $1 - C(u,v) \le (1-u) + (1-v) \le 2(1 - u \wedge v),$

for all $(u, v) \in [0, 1]^2$. Now take $(u, v), (s, t) \in \mathbb{E}$. If $s \le u$ and $t \le v$, then

$$|f_{(u,v)}(s,t)| = \frac{1 - C(u,v)}{q(u \wedge v)} \le 2(u \wedge v)^{-\omega} \le 2(s \wedge t)^{-\omega} = 2(s^{-\omega} \vee t^{-\omega}),$$

while if $s > u$ or $t > v$, then

$$|f_{(u,v)}(s,t)| = \frac{C(u,v)}{q(u \wedge v)} \le (1 - u \wedge v)^{-\omega}$$

$$\le (1 - s \vee t)^{-\omega} = (1-s)^{-\omega} \vee (1-t)^{-\omega}.$$

Hence, for every $f \in \mathcal{F}$ and every $(s,t) \in \mathbb{E}$, one has $|f(s,t)| \leq F(s,t)$. Furthermore, $P(F > x) \leq 4(x/2)^{-1/\omega} \wedge 1$ for every $x \geq 0$, so that the condition $\omega < 1/2$ ensures the integrability condition. Note in passing that since $Pf = 0$, one has $\mathbb{G}_n f = n^{1/2} \mathbb{P}_n f$ for all $f \in \mathcal{F}$. Hence, as $F$ is an envelope for $\mathcal{F}$, every sample path $f \mapsto \mathbb{G}_n f$ is an element of $\ell^\infty(\mathcal{F})$.

To prove (c), let $\mathbb{E}_0$ be a countable, dense subset of $\mathbb{E}$, and let $\mathcal{G}$ be the countable subset of $\mathcal{F}$ consisting of the zero function and the functions $f_{(u,v)}$ with $(u,v) \in \mathbb{E}_0$. Clearly, every $f \in \mathcal{F}$ is the pointwise limit of a sequence in $\mathcal{G}$. Furthermore, the envelope function $F$ is $P$-square integrable, and hence pointwise convergence in $\mathcal{F}$ implies $L^2(P)$ convergence. According to the definition at the bottom of page 116 of van der Vaart and Wellner (1996), this implies that $\mathcal{F}$ is indeed a pointwise separable class. $\square$

The limit process $\mathbb{G}$ whose existence is guaranteed by Lemma G.2 is a tight, Borel measurable element of $\ell^\infty(\mathcal{F})$ with Gaussian finite-dimensional distributions. To establish Theorem G.1, the idea is now to write $\mathbb{G}_{n,\omega} = T(\mathbb{G}_n)$ and $\mathbb{G}_\omega = T(\mathbb{G})$ as images of a continuous mapping $T \colon \ell^\infty(\mathcal{F}) \to \ell^\infty([0,1]^2)$ and to invoke the continuous mapping theorem.

To this end, introduce $\phi \colon [0,1]^2 \to \mathcal{F}$, which maps $(u,v) \in [0,1]^2$ to

$$\phi(u,v) = \begin{cases} f_{(u,v)}, & \text{if } u \wedge v \in (0,1), \\ 0, & \text{if } u = 0 \text{ or } v = 0 \text{ or } (u,v) = (1,1). \end{cases}$$

Now let $T \colon \ell^\infty(\mathcal{F}) \to \ell^\infty([0,1]^2)$ be defined by $T(z) = z \circ \phi$ for all $z \in \ell^\infty(\mathcal{F})$. One has $T(\mathbb{G}_n) = \mathbb{G}_{n,\omega}$ because if $(u,v) \in [0,1]^2$ with $u \wedge v \in (0,1)$, then

$$T(\mathbb{G}_n)(u,v) = \mathbb{G}_n f_{(u,v)} = n^{1/2}(\mathbb{P}_n - P)\frac{\mathbf{1}_{(0,u] \times (0,v]} - C(u,v)}{q(u \wedge v)} = \mathbb{G}_{n,\omega}(u,v),$$

while if $u = 0$, $v = 0$ or $(u,v) = (1,1)$ then

$$T(\mathbb{G}_n)(u,v) = \mathbb{G}_n 0 = 0 = \mathbb{G}_{n,\omega}(u,v).$$

The map $T$ is linear and bounded; it is thus continuous with respect to the topologies of uniform convergence on $\ell^\infty(\mathcal{F})$ and $\ell^\infty([0,1]^2)$. As for the map $\phi$, it is continuous with respect to the Euclidean metric on $[0,1]^2$ and the standard-deviation metric $\rho$ on $\mathcal{F}$, defined implicitly by

$$\rho^2(f,g) = \mathrm{E}\{(\mathbb{G}f - \mathbb{G}g)^2\} = Pf^2 - 2Pfg + Pg^2.$$

Indeed, the continuity of $\phi$ is easily derived from the bounds (G.1), together with the fact that for $f_{(u,v)}, f_{(s,t)}, 0 \in \mathcal{F}$,

$$\rho^2(f_{(u,v)}, 0) = \frac{C(u,v)\{1 - C(u,v)\}}{q(u \wedge v)^2}$$

and

$$\rho^2(f_{(u,v)}, f_{(s,t)}) = \frac{C(u,v)\{1 - C(u,v)\}}{q(u \wedge v)^2} + \frac{C(s,t)\{1 - C(s,t)\}}{q(s \wedge t)^2}$$

$$- 2\frac{C(u \wedge s, v \wedge t) - C(u,v)C(s,t)}{q(u \wedge v)q(s \wedge t)}$$

$$= \left(\frac{C(u,v)}{q(u \wedge v)} - \frac{C(u \wedge s, v \wedge t)}{q(s \wedge t)}\right)/q(u \wedge v)$$

$$+ \left(\frac{C(s,t)}{q(s \wedge t)} - \frac{C(u \wedge s, v \wedge t)}{q(u \wedge v)}\right)/q(s \wedge t)$$

$$- \left(\frac{C(u,v)}{q(u \wedge v)} - \frac{C(s,t)}{q(s \wedge t)}\right)^2.$$

The scene is finally set for an application of the continuous mapping theorem [see Theorem 1.3.6 in van der Vaart and Wellner (1996)].

PROOF OF THEOREM G.1. The continuous mapping theorem implies that $T(\mathbb{G}_n) \rightsquigarrow T(\mathbb{G})$ in $\ell^\infty([0,1]^2)$ as $n \to \infty$. The limit process $\mathbb{G}_\omega = T(\mathbb{G})$ is thus a tight, Borel measurable element of $\ell^\infty([0,1]^2)$ with the desired finite-dimensional distributions. Furthermore, $\phi$ is continuous, and according to the statements on page 41 of van der Vaart and Wellner (1996), the sample paths $f \mapsto \mathbb{G}f$ of $\mathbb{G}$ are almost surely uniformly continuous with respect to $\rho$. Therefore, the sample paths $(u,v) \mapsto \mathbb{G}_\omega(u,v) = \mathbb{G}\phi(u,v)$ are almost surely continuous, as well. Hence, the process $\mathbb{G}_\omega$ has all the stated properties. $\square$

Département de mathématiques
et de statistique
Université Laval
1045, avenue de la Médecine
Québec, Canada G1V 0A6
E-mail: Christian.Genest@mat.ulaval.ca

Institut de statistique
Université catholique de Louvain
Voie du Roman Pays, 20
B-1348 Louvain-la-Neuve
Belgium
E-mail: Johan.Segers@uclouvain.be